\documentclass[aop,MSNbibl,citesort,seceqn,dvips]{arximspdf}

\doi{10.1214/11-AOP683}
\volume{40}
\issue{6}
\pubyear{2012}
\firstpage{2400}
\lastpage{2438}

\makeatletter

\newcommand{\pmb}{\mathbf}
\newcommand{\ulN}{\displaystyle\mathop{\longrightarrow}_{N}}

\newcommand{\IC}{{\mathbb C}}
\newcommand{\IE}{{\mathbb E}}
\newcommand{\IP}{{\mathbb P}}
\newcommand{\IR}{{\mathbb R}}
\newcommand{\IN}{{\mathbb N}}
\newcommand{\IZ}{{\mathbb Z}}

\newcommand{\fr}{\frac{1}{2}}
\newcommand{\point}{{\bolds{.}}}
\newcommand{\wh}{\widehat}
\newcommand{\wt}{\widetilde}
\newcommand{\cA}{\mathcal{A}}
\newcommand{\cE}{\mathcal{E}}
\newcommand{\cL}{\mathcal{L}}
\newcommand{\cI}{\mathcal{I}}
\newcommand{\cS}{\mathcal{S}}

\newtheorem{theorem}{Theorem}[section]
\newtheorem{lemma}[theorem]{Lemma}

\newproclaim{remark}[theorem]{Remark}

\makeatother

\begin{document}
\begin{frontmatter}

\title{Random interlacements and the Gaussian free field}
\runtitle{Random interlacements and the Gaussian free field}

\begin{aug}
\author[A]{\fnms{Alain-Sol} \snm{Sznitman}\corref{}\thanksref{t1}\ead[label=e1]{sznitman@math.ethz.ch}}
\runauthor{A.-S. Sznitman}
\affiliation{ETH Zurich}
\address[A]{Departement Mathematik \\
ETH-Zurich\\
CH-8092 Z\"urich\\
Switzerland\\
\printead{e1}} 
\end{aug}

\thankstext{t1}{Supported in part by the Grant ERC-2009-AdG 245728-RWPERCRI.}

\received{\smonth{2} \syear{2011}}
\revised{\smonth{5} \syear{2011}}

%
\begin{abstract}
We consider continuous time random interlacements on $\IZ^d$, $d \ge
3$, and characterize the distribution of the corresponding stationary
random field of occupation times. When $d=3$, we relate this random
field to the two-dimensional Gaussian free field pinned at the origin
by looking at scaled differences of occupation times of long rods by
random interlacements at appropriately tuned levels. In the main
asymptotic regime, a scaling factor appears in the limit, which is
independent of the free field, and distributed as the time-marginal of
a zero-dimensional Bessel process. For arbitrary $d \ge3$, we also
relate the field of occupation times at a level tending to infinity, to
the $d$-dimensional Gaussian free field.
\end{abstract}

%

\end{frontmatter}

\setcounter{section}{-1}
\section{Introduction}\label{intro}

In this article we consider continuous time random interlacements on
$\IZ^d$, $d \ge3$, where each doubly infinite trajectory modulo
time-shift in the interlacement is decorated by i.i.d. exponential
variables with parameter $1$ which specify the time spent by the
trajectory at each step. We are interested in the random field of
occupation times, that is, the total time spent at each site of $\IZ
^d$ by the collection of trajectories with label at most $u$ in the
interlacement point process.

When $d = 3$, we relate this stationary random field to the
\textit{two-dimensional} Gaussian free field pinned at the origin by looking at
the properly scaled field of differences of occupation times of long
rods of size $N$, when the level $u$ is either proportional to $\log
N/N$ or much larger than $\log N/N$. The choice of $u$ proportional to
$\log N/N$ corresponds to a nondegenerate probability that the
interlacement at level $u$ meets a given rod. In the asymptotic regime
it brings into play an independent proportionality factor of the
Gaussian free field, which is distributed as a certain time-marginal of
a zero-dimensional Bessel process. This random factor disappears from
the description of the limiting random field, when instead $u N / \log
N$ tends to infinity.

For arbitrary $d \ge3$, we also relate the properly scaled field of
differences of occupation times of sites by the interlacement at a
level $u$ tending to infinity, with the Gaussian free field on $\IZ^d$.

Rather than discussing our results any further, we first present the
model and refer to Section~\ref{sec1} for additional details. We consider the
spaces $\wh{W}_+$ and $\wh{W}$ of infinite and doubly infinite $\IZ
^d \times(0,\infty)$-valued sequences, with $d \ge3$, such that the
$\IZ^d$-valued components form an infinite, respectively, doubly
infinite, nearest neighbor trajectory spending finite time in any
finite subset of $\IZ^d$, and such that the $(0,\infty)$-valued
components have an infinite sum in the case of $\wh{W}_+$, and
infinite ``forward'' and ``backward'' sums, when restricted to positive
and negative indices, in the case of $\wh{W}$.

We write\vspace*{1pt} $X_n, \sigma_n$, with $n \ge0$, or $n \in\IZ$, for the
respective $\IZ^d$- and $(0,\infty)$-valued canonical coordinates on
$\wh{W}_+$ and $\wh{W}$. We denote by $P_x$, $x \in\IZ^d$, the law
on $\wh{W}_+$ endowed with the canonical $\sigma$-algebra, under
which $X_n$, $n \ge0$, are distributed as simple random walk starting
at $x$, and $\sigma_n$, $n \ge0$, are i.i.d.\vspace*{1pt} exponential variables
with parameter $1$, independent from the $X_n$, $n \ge0$. We write
$\wh{W}^*$ for the space $\wh{W}$ modulo time-shift, that is, $\wh
{W}^* = \wh{W}/\sim$, where for $\wh{w}, \wh{w}'$ in $\wh{W}$,
$\wh{w} \sim\wh{w}'$ means that $\wh{w}(\cdot) = \wh{w}'(\cdot+
k)$ for some $k \in\IZ$. We denote\vspace*{1pt} by $\pi^*$: $\wh{W} \rightarrow\wh
{W}^*$ the canonical map, and endow $\wh{W}^*$ with the $\sigma
$-algebra consisting\vspace*{2pt} of sets with an inverse image under $\pi^*$
belonging to the canonical $\sigma$-algebra of $\wh{W}$.

The continuous time interlacement point process on $\IZ^d$, $d \ge3$,
is a Poisson point process on $\wh{W} \times\IR_+$. Its intensity
measure has the form $\wh{\nu}(d \wh{w}^*) \,du$, where $\wh{\nu}$
is the $\sigma$-finite measure on $\wh{W}^*$ such that for any finite
subset $K$ of $\IZ^d$, the restriction of $\wh{\nu}$ to the subset
of $\wh{W}^*$ made of $\wh{w}^*$ for which the $\IZ^d$-valued
trajectory modulo time-shift enters $K$, is equal to $\pi^* \circ\wh
{Q}_K$, the image of $\wh{Q}_K$ under $\pi^*$, where $\wh{Q}_K$ is
the finite measure specified by:
\begin{eqnarray}\label{0.1}
&&\begin{tabular}{p{320pt}}
\mbox{}\quad \hphantom{i}(i) $\wh{Q}_K(X_0 = x) = e_K(x)$, with $e_K$ the
equilibrium measure of $K$; see~(\ref{1.4});
\end{tabular}\hspace*{-28pt}
\nonumber\\[-8pt]\\[-8pt]
&&\begin{tabular}{p{320pt}}
\mbox{}\quad (ii) when $e_K(x) > 0$, conditionally on $X_0 = x$, $(X_n)_{n \ge0}$,
$(X_{-n})_{n \ge0}$, $(\sigma_n)_{n \in\IZ}$ are independent,
respectively, distributed as simple random walk starting at $x$, as
simple random walk starting at $x$ conditioned never to return to $K$
and as a doubly infinite sequence of independent exponential
variables with parameter $1$.
\end{tabular}\hspace*{-28pt}
\nonumber
\end{eqnarray}
The existence and uniqueness of such a measure $\wh{\nu}$ can be shown
just as in Section 1 of~\cite{Szni10a}. The canonical continuous time
interlacement point process is then constructed on a space $(\Omega,
\cA, \IP)$, similar to (1.16)\vspace*{1pt} of~\cite{Szni10a}, with $\omega= \sum_{i
\ge0} \delta_{(\wh{w}_i^*, u_i)}$ denoting a generic element of the
set $\Omega$. We also refer to Remark~\ref{rem2.4}(4) which explains
how $\IZ^d$, $d \ge3$, can be replaced with a transient weighted
graph, and continuous time random interlacements on a transient
weighted graph are constructed.

In the present work our main interest focuses on the collection of
(continuous) occupation times
%
%
\begin{eqnarray}\label{0.2}
&&L_{x,u}(\omega) = \sum_{i \ge0} \sum_{n \in\IZ} \sigma_n(\wh
{w}_i) 1\{X_n(\wh{w}_i) = x, u_i \le u\}\qquad\mbox{for } x \in
\IZ^d, u \ge0
\nonumber\\[-8pt]\\[-8pt]
&&\eqntext{\mbox{where }\displaystyle  \omega= \sum_{i \ge0} \delta_{(\wh{w}^*_i,u_i)} \in
\Omega \mbox{ and } \pi^*(\wh{w}_i) = \wh{w}_i^* \mbox{ for
each $i \ge0$}.}
\end{eqnarray}
We compute the Laplace functional of this random field and show in
Theorem~\ref{theo2.3} that when $V$ is a nonnegative function on $\IZ
^d$ with finite support, one has the identity
%
%
\begin{equation}\label{0.3}\qquad
\IE\biggl[\exp\biggl\{ - \sum_{x \in\IZ^d} V(x) L_{x,u}\biggr\}\biggr]
= \exp\biggl\{ - u \frac{\sum_{\phi\not= I} c_I \Pi_{x
\in I} V(x)}{\sum_I g_I \Pi_{x \in I} V(x)}\biggr\},\qquad u \ge0,
\end{equation}
where, in the above formula, $I$ runs over the collection of subsets of
the support of~$V$, $g_I$ denotes the determinant of the Green function
$g(\cdot,\cdot)$ restricted to $I \times I$ [see~(\ref{1.1})] and
$c_I$ the sum of the coefficients of the matrix of cofactors of the
above matrix. Both quantities are positive and [see~(\ref{2.13})] their
ratio $c_I/g_I$ coincides with the capacity of $I$, that is, the total
mass of the equilibrium measure $e_I$ of $I$. We refer to~(\ref{2.26})
for the extension of this formula to the case where $\IZ^d$ is
replaced by a transient weighted graph. One can also consider the
discrete occupation times, where $\sigma_n$ is replaced by $1$ in~(\ref{0.2}); however, this random field turns out to be somewhat less
convenient to handle than $(L_{x,u})_{x \in\IZ^d}$ for the kind of
questions we investigate here; see Remark~\ref{rem2.4}(5).

The continuous time interlacement point process is related to the
Poisson point process of Markov loops initiated in~\cite{Syma69},
which later found various incarnations (see, e.g., Theorem 2.1
of~\cite{BrydFrohSpen82}, Sections 4 and 3 of~\cite{Dynk83}, and \cite
{LawlWern04}, Chapter 9 of~\cite{LawlLimi10}) and was extensively
analyzed in~\cite{Leja10,Leja11}. Heuristically random
interlacements correspond to a ``restriction to loops going through
infinity'' of this Poisson point process; see~\cite{Leja11}, page 85.
It has been shown in Theorem 13 of~\cite{Leja10} (see also \cite
{Leja11}, page 61) that the field of occupation times of the Poisson
point process of Markov loops on a finite weighted graph with
nondegenerate killing, at a suitable choice of the level is
distributed as half the square of the Gaussian free field on the finite
graph. No such identity holds in our context when considering a fixed
level $u$; see Remark~\ref{rem2.4}. However, and this is the main
object of this article, we present limiting procedures which relate the
field of occupation times of random interlacements to the Gaussian free field.

The link with the two-dimensional Gaussian free field comes as follows.
We look at the occupation times of long vertical rods in $\IZ^3$, by
random interlacements at properly tuned levels. Let us incidentally
mention that the consideration of long rods in the context of random
interlacements has been helpful in several instances, for example,
Section 3 of~\cite{SidoSzni09a} or Section~5 of~\cite{Szni10b}. They
typically have been used as a tool in the detection of long
$*$-crossings in planes, left by the trajectories of the random
interlacements at level $u$, and have enabled us to quantify the rarity
of such crossings when $u$ is small. Here the rods in question are the
subsets of $\IZ^3$,
%
%
\begin{equation}\label{0.4}
J_y = \{x = (y,k) \in\IZ^3; 1 \le k \le N\}\qquad \mbox{for $y \in
\IZ^2$ and $N>1$},
\end{equation}
and the corresponding $\IZ^2$-stationary field of occupation times is
given by
%
%
\begin{equation}\label{0.5}
\cL_{y,u} = \sum_{x \in J_y} L_{x,u},\qquad y \in\IZ^2, u \ge0.
\end{equation}
We choose the levels $(u_N)_{N > 1}$ and $(u'_N)_{N>1}$, so that
%
%
\begin{equation}\label{0.6}
\mbox{(i)\quad} u_N = \alpha\frac{\log N}{N} \qquad \mbox{with
$\alpha> 0$},\qquad
\mbox{(ii)\quad} \frac{\log N}{N} = o(u'_N).
\end{equation}
The choice in~(\ref{0.6})(i) corresponds to a nondegenerate limiting
probability $\exp\{-\frac{\pi}{3} \alpha\}$ that the
interlacement at level $u_N$ does not meet any given rod $J_y$ [see
(\ref{4.74})] whereas the choice in~(\ref{0.6})(ii) induces a
vanishing limit for the corresponding probability.

If we now introduce the Gaussian free field pinned at the origin, or
more precisely [see~(\ref{1.29})] a centered Gaussian field $(\psi
_y)_{y \in\IZ^2}$, with covariance $3(a(y) + a(y') - a(y' - y))$,
$y,y' \in\IZ^2$, where $a(\cdot)$ is the potential kernel of the
two-dimensional simple random walk [see~(\ref{1.6})] and $R$ an
independent nonnegative random variable, having the law $\mathrm{BES}^o(\sqrt
{\alpha}, \frac{3}{2 \pi})$ of a zero-dimensional Bessel process at
time $\frac{3}{2 \pi}$ starting in $\sqrt{\alpha}$ at time $0$ [see
(\ref{1.30})]\vspace*{1pt} we show in Theorems~\ref{theo4.2} and~\ref{theo4.9}
that when $N$ tends to infinity,
%
%
\begin{equation} \label{0.7}
\biggl(\frac{\cL_{y,u_N}}{\log N}\biggr)_{y \in\IZ^2} \mbox
{ converges in distribution to the flat field with value $R^2$}\hspace*{-35pt}
\end{equation}
and that
\begin{eqnarray} \label{0.8}
&&\biggl(\frac{\cL_{y,u_N}- \cL_{0,u_N}}{\sqrt{\log N}}\biggr)_{y
\in\IZ^2} \mbox{ converges in distribution to the random}\nonumber\\[-8pt]\\[-8pt]
&&\mbox{field
$(R\psi_y)_{y \in\IZ^2}$.}\nonumber
\end{eqnarray}
In the case~(\ref{0.6}) (ii) we instead find that when $N$ goes to infinity,
%
%
\begin{equation} \label{0.9}
\biggl(\frac{\cL_{y,u'_N}}{Nu'_N}\biggr)_{y \in\IZ^2} \mbox
{ converges in distribution to the flat field with value $1$}\hspace*{-35pt}
\end{equation}
and that
%
%
\begin{equation} \label{0.10}
\biggl(\frac{\cL_{y,u'_N}- \cL_{0,u'_N}}{\sqrt{N u'_N}}
\biggr)_{y \in\IZ^2} \mbox{ converges in distribution to $(\psi_y)_{y
\in\IZ^2}$}.
\end{equation}
There is an important connection between random interlacements and the
structure left locally by a random walk on a large torus; see
\cite{Wind08,TeixWind10}. In this light one may wonder whether
some of the above results have counterparts in the case of a simple
random walk on a large two-dimensional torus. We refer to Remark
\ref{rem4.10}(1) for more on this issue. Some consequences of the
above limit results for discrete occupation times of long rods can also
be found in Remark~\ref{rem4.10}(2).

In this article, we provide yet a further link between random
interlacements and the Gaussian free field, by considering the
occupation times of random interlacements at a level $u$ tending to
infinity. If $(\gamma_x)_{x \in\IZ^d}$ stands for the Gaussian free
field on $\IZ^d$, $d \ge3$, that is, the centered Gaussian field with
covariance function $E[\gamma_x \gamma_{x'}] = g(x,x')$, $x,x' \in
\IZ^d$, we show in Theorem~\ref{theo5.1} that when $u$ tends to infinity,
\begin{eqnarray}\label{0.11}
&&\biggl(\frac{1}{u} L_{x,u}\biggr)_{x \in\IZ^d} \mbox
{converges in distribution toward the flat}\nonumber\\[-8pt]\\[-8pt]
&&\mbox{field with value $1$}\nonumber
\end{eqnarray}
and that
%
%
\begin{equation} \label{0.12}
\biggl(\frac{L_{x,u} - L_{x,0}}{\sqrt{2 u}}\biggr)_{x \in\IZ^d}
\mbox{ converges in distribution toward $(\gamma_x - \gamma_0)_{x
\in\IZ^d}$}.\hspace*{-35pt}
\end{equation}
We refer to Remark~\ref{rem5.2} for the extension of these results to
the case of random interlacements on a transient weighted graph and to
discrete occupation times.

Let us say a few words concerning proofs. We provide in Theorem \ref
{theo2.1} an expression for the characteristic function of $\sum_{x
\in\IZ^d} V(x) L_{x,u}$, with $V$ finitely supported, which shows
that close to the origin it can be expressed as the exponential of an
analytic function. This identity on the one hand leads to~(\ref{0.3});
see Theorem~\ref{theo2.3}. On the other hand, this identity underlies
the general line of attack, which we employ when proving the limit
theorems corresponding to~(\ref{0.7})--(\ref{0.10}) and (\ref
{0.11}),~(\ref{0.12}). Namely we investigate the asymptotic behavior
of the power series representing the above mentioned analytic
functions. The proof of~(\ref{0.8}) is by far the most delicate. We
analyze the large $N$ asymptotics of the power series, expressing the
logarithm of the characteristic function of $\sum_{y \in\IZ^2} W(y)
\cL_{y,u_N}$ close to the origin, with $W$ finitely supported on $\IZ
^2$, and such that $\sum_y W(y) = 0$. This asymptotic analysis relies
on certain cancellations, which take place and enable us to control the
coefficients of the power series. In the crucial Theorem~\ref{theo4.1}
we bound these coefficients, show the asymptotic vanishing of odd
coefficients and compute the (nonvanishing) limit of even
coefficients. This theorem contains enough information to yield both
(\ref{0.8}) and~(\ref{0.10}); see Theorem~\ref{theo4.2}. Once (\ref
{0.8}),~(\ref{0.10}) are proved,~(\ref{0.7}),~(\ref{0.9}) follow in
a simpler fashion and in essence only require the consideration of one
single rod, say $J_0$. The proof of~(\ref{0.11}),~(\ref{0.12}) in
Theorem~\ref{theo5.1} follows a similar pattern, but is substantially simpler.

Let us now describe how this article is organized.

In Section~\ref{sec1} we provide additional notation and collect some results
concerning potential theory, the two-dimensional free field and
zero-dimensional Bessel processes.

Section~\ref{sec2} contains the identity for the characteristic functional of
the field of occupation times in Theorem~\ref{theo2.1} and the proof
of formula~(\ref{0.3}) for the Laplace functional in Theorem \ref
{theo2.3}. The extension of these results to the set-up of weighted
graphs can be found in Remark~\ref{rem2.4}(4).

In Section~\ref{sec3} we collect estimates as preparation for the study in the
next section of occupation times of long rods in $\IZ^3$.

Section~\ref{sec4} presents the limiting results~(\ref{0.7})--(\ref{0.10})
(see Theorems~\ref{theo4.2} and~\ref{theo4.9}) relating random
interlacements in $\IZ^3$ to the two-dimensional free field. The heart
of the matter lies in Theorem~\ref{theo4.1}, where controls over the
relevant power series are derived.

In Section~\ref{sec5} we prove~(\ref{0.11}),~(\ref{0.12}) in Theorem \ref
{theo5.1} and provide in Remark~\ref{rem5.2} the extension of these
results to the case of transient weighted graphs, and to discrete
occupation times.

Finally let us explain our convention concerning constants. We denote
with $c,c',\wt{c},\overline{c}$ positive constants changing from
place to place. Numbered constants refer to the value corresponding to
their first appearance in the text. In Sections~\ref{sec1},~\ref{sec2} and
\ref{sec5} constants
only depend on $d$. In Section~\ref{sec3}, where $d=3$, they depend on $\Lambda
$ in~(\ref{3.1}), and in Section~\ref{sec3}, where $d=3$ as well, on $\Lambda$
and $W$; see~(\ref{4.3}). Otherwise dependence of constants on
additional parameters appears in the notation.

\section{Notation and some useful facts}\label{sec1}

In this section we provide some additional notation and recall various
useful facts concerning random walks, discrete potential theory, the
two-dimensional free field and zero-dimensional Bessel processes.

We let $\IN= \{0,1,\ldots\}$ denote the set of natural numbers. When
$u$ is a nonnegative real number we let $[u]$ stand for the integer
part of $u$. Given a finite set $A$, we denote by $|A|$ its
cardinality. We write \mbox{$|\cdot|$} for the Euclidean norm on $\IR^d$, $d
\ge1$. For $A, A^\prime\subseteq\IZ^d$, we denote by $d(A,A^\prime
) = \inf\{|x-x^\prime|$; $x \in A$, $x^\prime\in A^\prime\}$ the
mutual distance of $A$ and $A^\prime$. When $A = \{x\}$, we write
$d(x,A^\prime)$ in place of $d(A,A^\prime)$ for simplicity. We write
$U \subset\subset\IZ^d$, to indicate that $U$ is a finite subset of
$\IZ^d$. Given $f,g$ square summable functions on $\IZ^d$ we write
$(f,g) = \sum_{x \in\IZ^d} f(x) g(x)$ for their scalar product. When
$U \subseteq\IZ^d$, and $f$ is a function on $U$, we routinely
identify $f$ with the function on $\IZ^d$, which vanishes outside $U$
and coincides with $f$ on $U$. We denote the sup-norm of such a
function with $\|f\|_{L^\infty(U)}$, and sometimes with $\|f\|_\infty
$, when there is no ambiguity.

Given $U \subseteq\IZ^d$, we write $H_U = \inf\{n \ge0; X_n \in U\}
$, $\wt{H}_U = \inf\{n \ge1; X_n \in U\}$ and $T_U = \inf\{n \ge0;
X_n \notin U\}$ for the entrance time of $U$, the hitting time of $U$,
and the exit time from $U$. When $\rho$ is a measure on $\IZ^d$, we
denote by $P_\rho$ the measure $\sum_{x \in\IZ^d} \rho(x) P_x$,
and by $E_\rho$ the corresponding expectation. So far $P_x$, $x \in
\IZ^d$, has only been defined when $d \ge3$; see above~(\ref{0.1}).
When $d = 1$ or $2$, this notation simply stands for the canonical law
of simple random walk starting at $x$, and $X_n$, $n \ge0$, for the
canonical process.

When $d \ge3$, we denote by $g(\cdot,\cdot)$ the Green function
%
%
\begin{equation}\label{1.1}
g(x,x^\prime) = \sum_{n \ge0} P_x [X_n = x^\prime]\qquad \mbox{for
$x,x^\prime$ in $\IZ^d$}.
\end{equation}
It is a symmetric function, and due to translation invariance one has
%
%
\begin{equation}\label{1.2}
g(x,x^\prime) = g(x^\prime-x ) = g(x-x^\prime)\qquad\mbox{where
$g(\cdot) = g(\cdot,0)$}.
\end{equation}
Classically one knows that $g(\cdot) \le g(0)$, and that (see
\cite{Lawl91}, page 31)
%
%
\begin{equation}\label{1.3}
c^\prime(1 \vee|x|)^{-(d-2)} \le g(x) \le c(1 \vee|x|)^{-(d-2)}\qquad
\mbox{for $x \in\IZ^d$}.
\end{equation}
When $K \subset\subset\IZ^d$, the equilibrium measure of $K$ and the
capacity of $K$, that is, the total mass of $e_K$, are denoted by
\begin{eqnarray}\label{1.4}
e_K(x) &=& P_x[\wt{H}_K = \infty] 1_K (x)\qquad \mbox{for $x \in\IZ
^d$}\quad \mbox{and}\nonumber\\[-8pt]\\[-8pt]
\operatorname{cap}(K) &=& \sum_{x \in\IZ^d} P_x[\wt{H}_K =
\infty].\nonumber
\end{eqnarray}
One can express the probability to enter $K$ via the formula
%
%
\begin{equation}\label{1.5}
P_x[H_K < \infty] = \sum_{x^\prime\in K} g(x,x^\prime)
e_K(x^\prime)\qquad \mbox{for $x \in\IZ^d$}.
\end{equation}
We will also consider the two-dimensional potential kernel (see (1.40),
page~37 of~\cite{Lawl91}, or pages 121, 122, 148 of~\cite{Spit01})
%
%
\begin{equation}\label{1.6}
a(y) = \lim_{n \rightarrow\infty} \sum^n_{j=0} P_0 [X_j = 0] -
P_0[X_j = y]\qquad \mbox{for $y \in\IZ^2$} .
\end{equation}
It is a nonnegative function on $\IZ^2$, which is symmetric and
satisfies (cf. Proposition~P2, page 123 of~\cite{Spit01})
%
%
\begin{equation}\label{1.7}
\lim_{y^\prime\rightarrow\infty} a(y + y^\prime) - a(y^\prime) =
0\qquad \mbox{for any $y \in\IZ^2$}.
\end{equation}
In Sections~\ref{sec3} and~\ref{sec4} [see also~(\ref{0.4})] we consider long vertical
rods, which are the subsets of $\IZ^3$ defined for $y \in\IZ^2$ and
$N > 1$, by
%
%
\begin{equation}\label{1.8}
J_y = \{y\} \times J \subseteq\IZ^3\qquad \mbox{where $J = \{1,\ldots, N\}$}.
\end{equation}
The next lemma collects limit statements concerning the potentials of
long rods, and in particular relates the difference of such potentials
to the two-dimensional potential kernel.
\begin{lemma}[$(d=3, N > 1, y \in\IZ^2)$]\label{lem1.1}
%
%
\begin{equation}\label{1.9}\quad
\lim_N \frac{1}{2 \log N} \sum_{|z| \le N}
g((0,z)) =
\frac{3}{2 \pi}\qquad \mbox{[with $z \in\IZ$ and $(0,z) \in\IZ^3$]}.
\end{equation}
For $x = (0,z)$ in $J_0$ and $y \in\IZ^2$, one has
%
%
\begin{equation}\label{1.10}
\sum_{x^\prime\in J_0} g(x,x^\prime) - \sum_{x^{\prime\prime} \in
J_y} g(x,x^{\prime\prime}) =
\frac{3}{2} a(y) - b_N(y,z) ,
\end{equation}
where $b_N$ is a nonnegative function on $\IZ^2 \times J$ such that
%
%
\begin{eqnarray}\label{1.11}
b_N(y,z) \le\psi_y(d(z,J^c))\nonumber\\[-8pt]\\[-8pt]
&&\eqntext{\mbox{where }\displaystyle \lim_{r \rightarrow
\infty} \psi_y(r) = 0\qquad \mbox{for each $y \in\IZ^2$}.}
\end{eqnarray}
\end{lemma}
\begin{pf}
Claim~(\ref{1.9}) is an immediate consequence of the fact that $\sum
^N_1 \frac{1}{k} \sim\log N$, as $N$ goes to infinity and
(cf. Theorem 1.5.4, page 31 of~\cite{Lawl91})
%
%
\begin{equation}\label{1.12}
g(x) \sim\frac{3}{2 \pi} |x|^{-1}\qquad \mbox{as } x \rightarrow\infty.
\end{equation}
We now turn to the proof of~(\ref{1.10}),~(\ref{1.11}). We denote by
$\wt{Y}_\point$ and $\wt{Z}_\point$ independent continuous time
random walks on $\IZ^2$ and $\IZ$ with respective\vspace*{1pt} jump rates 2 and 1,
starting at the respective origins of $\IZ^2$ and $\IZ$. So $(\wt
{Y}_\point, \wt{Z}_\point)$ is a continuous time random walk on $\IZ
^3$, starting at the origin, with jump rate equal to 3, and the
left-hand side of~(\ref{1.10}) equals
%
%
\begin{eqnarray}\label{1.13}
&&
3 E\biggl[\int_0^\infty1\{\wt{Y}_s = 0, \wt{Z}_s + z \in J\}
\,ds - \int^\infty_0 1\{\wt{Y}_s = y, \wt{Z}_s + z \in J\}
\,ds\biggr] \nonumber\\
&&\qquad\stackrel{\mathrm{independence}}{=}
3 \int^\infty_0 (P[\wt{Y}_s = 0] - P[\wt{Y}_s = y]) P[\wt
{Z}_s + z \in J] \,ds \\
&&\hspace*{19pt}\qquad= 3(I_1 - I_2),
\nonumber
\end{eqnarray}
where we have set
\begin{eqnarray}\label{1.14}
I_1 & = & \int^\infty_0 P[\wt{Y}_s = 0] - P[\wt{Y}_s = y] \,ds,
\nonumber\\[-8pt]\\[-8pt]
I_2 & = & \int^\infty_0 (P[\wt{Y}_s = 0] - P[\wt{Y}_s = y])
P[\wt{Z}_s + z \notin J] \,ds,
\nonumber
\end{eqnarray}
and we note that the integrand in $I_1$ is nonnegative as a direct
application of the Chapman--Kolmogorov equation at time $\frac{s}{2}$
and the Cauchy--Schwarz inequality. If we let $Y_k, k \ge0$ and $T_k,
k \ge0$, (with $T_0 = 0$), stand for the discrete skeleton of $\wt
{Y}_\point$ and its successive jump times, we see that for $T > 0$,
\begin{eqnarray}\label{1.15}
&&\int^T_0 P[\wt{Y}_s = 0] - P[\wt{Y}_s = y] \,ds \nonumber\\
&&\hspace*{18.68pt}\qquad= \sum_{k \ge
0} E[(T_{k+1} \wedge T - T_k \wedge T) 1\{Y_k = 0\}]
\nonumber\\[-8pt]\\[-8pt]
&&\hspace*{18.68pt}\qquad\quad{}-E[(T_{k+1} \wedge T - T_k \wedge T) 1\{Y_k = y\}]\nonumber\\
&&\qquad\stackrel{\mathrm{independence}}{=}
\sum_{k \ge0} E[T_{k+1} \wedge T - T_k \wedge T] (P[Y_k = 0] - P[Y_k
= y] ).
\nonumber
\end{eqnarray}
Observe that $T_{k+1} - T_k$ is an exponential variable with parameter
$2$, which is independent from $T_k$, so that for $k \ge0$
%
%
\begin{equation}\label{1.16}\qquad\quad
a_{k,T} \stackrel{\mathrm{def}}{=} E[T_{k+1} \wedge T - T_k \wedge T] =
E\biggl[T_k \le T, 2 \int^\infty_0 s \wedge(T-T_k) e^{-2s}
\,ds\biggr]
\end{equation}
decreases to zero as $k$ tends to infinity and increases to $\frac
{1}{2}$ as $T$ tends to $\infty$. We set $s_k = \sum_{0 \le j \le k}
P[Y_j = 0] - P[Y_j = y]$, for $k \ge0$, so that by~(\ref{1.6}), $\lim
_k s_k = a(y)$. After summation by parts in the last member of (\ref
{1.15}), we find that
%
%
\begin{equation}\label{1.17}
\int^T_0 P[\wt{Y}_s = 0] - P[\wt{Y}_s= y] \,ds = \sum_{k \ge
0} (a_{k,T} - a_{k+1,T}) s_k .
\end{equation}
Using the observations below~(\ref{1.16}) we see that the left-hand
side of~(\ref{1.17}) tends to $\frac{1}{2} a(y)$ as $T$ goes to
infinity, so that
%
%
\begin{equation}\label{1.18}
I_1 = \tfrac{1}{2} a(y) .
\end{equation}
As for $I_2$, which is nonnegative due to the remark below (\ref
{1.14}), we see that $3 I_2 \le\psi_y (d(z,J^c))$, where we have set
%
%
\begin{equation}\label{1.19}\qquad
\psi_y(r) = 3 \int^\infty_0 (P[\wt{Y}_s = 0] - P[\wt{Y}_s =
y]) P[ | \wt{Z}_s | \ge r] \,ds\qquad \mbox{for $r \ge0$}.
\end{equation}
If is plain that $\psi_y$ is a nonincreasing function, which tends to
zero at infinity by dominated convergence. This completes the proof of
Lemma~\ref{lem1.1}.
\end{pf}

We now turn to the discussion of the two-dimensional massless Gaussian
free field pinned at the origin. For this purpose we begin by the
consideration of the more traditional two-dimensional massless Gaussian
free field with Dirichlet boundary conditions outside the square $U_L =
[-L,L]^2$, with $L \ge1$; see, for instance,~\cite{BoltDeusGiac01}. It
is a centered Gaussian field $\varphi_{y,L}$, $y \in\IZ^2$, with
covariance function
%
%
\begin{equation}\label{1.20}
E[\varphi_{y,L} \varphi_{y^\prime,L}] = g_L(y,y^\prime)\qquad
\mbox{for $y,y^\prime\in\IZ^2$},
\end{equation}
where $g_L(\cdot,\cdot)$ stands for the Green function of the
two-dimensional random walk killed when exiting $U_L$
%
%
\begin{equation}\label{1.21}
g_L(y,y^\prime) = E_y \biggl[\sum_{k \ge0} 1 \{X_k = y^\prime, k <
T_{U_L}\}\biggr]\qquad \mbox{for } y,y^\prime\in\IZ^2 .
\end{equation}
Writing $H_0$ in place of $H_{\{0\}}$ [see above~(\ref{1.1})] it
follows from the strong Markov property and~(\ref{1.20}), (\ref
{1.21}), that for any $y \in\IZ^2$,
\[
\varphi_{y,L} - P_y[H_0 < T_{U_L}] \varphi_{0,L} \mbox{ is
orthogonal to $\varphi_{0,L}$}.
\]
Hence defining for any $\gamma\in\IR$,
%
%
\begin{equation}\label{1.22}\qquad
\Phi_{y,L}(\gamma) = \varphi_{y,L} - P_y[H_0 < T_{U_L}] \varphi
_{0,L} + P_y [H_0 < T_{U_L}] \gamma,\qquad y\in\IZ^2 ,
\end{equation}
the law of the above random field is a regular conditional probability
for the law of $(\varphi_{y,L})$ given its value at the origin
$\varphi_{0,L} = \gamma$. The next lemma will provide two possible
interpretations for the centered Gaussian field we consider in the
sequel, in terms of the two-dimensional massless Gaussian free field.
\begin{lemma}\label{lem1.2}
For $y,y^\prime$ in $\IZ^2$, one has
\begin{eqnarray}\label{1.23}
a(y) + a(y^\prime) - a(y^\prime- y) & = & \lim_{L \rightarrow\infty}
E[(\varphi_{y,L} - \varphi_{0,L}) (\varphi_{y^\prime,L} - \varphi_{0,L})]
\nonumber\\[-8pt]\\[-8pt]
& = & \lim_{L \rightarrow\infty} E[\Phi_{y,L} (0) \Phi_{y^\prime
,L}(0)] .
\nonumber
\end{eqnarray}
\end{lemma}
\begin{pf}
By~(\ref{1.20}) we see that
\begin{eqnarray}\label{1.24}
&&
E[(\varphi_{y,L} - \varphi_{0,L}) (\varphi_{y',L} - \varphi_{0,L})]
\nonumber\\
&&\qquad= g_L(y,y^\prime) + g_L(0,0) - g_L(0,y) - g_L(0,y^\prime)
\nonumber\\[-8pt]\\[-8pt]
&&\qquad=
g_L(y,y) - g_L(0,y) + g_L(0,0) - g_L(0,y^\prime) \nonumber\\
&&\qquad\quad{} - \bigl(g_L(y,y) -
g_L(y,y^\prime)\bigr) .
\nonumber
\end{eqnarray}
From Proposition 1.6.3, page 39 of~\cite{Lawl91}, one knows that for
$y_1,y_2$ in $U_L$,
%
%
\begin{equation}\label{1.25}
g_L(y_1,y_2) = \sum_{y^\prime\in\partial U_L} P_{y_1} [X_{T_{U_L}} =
y^\prime] a(y^\prime- y_2) - a(y_1 - y_2),\vadjust{\goodbreak}
\end{equation}
so that by~(\ref{1.7}) and $a(0) = 0$, we find that
%
%
\begin{equation}\label{1.26}
\lim_{L \rightarrow\infty} g_L(y_1,y_1) - g_L(y_1,y_2) = a(y_1 -
y_2)\qquad \mbox{for } y_1,y_2 \in\IZ^2.
\end{equation}
Coming back to~(\ref{1.24}) and keeping in mind the symmetry of
$g_L(\cdot,\cdot)$, the first equality of~(\ref{1.23}) follows. As
for the second equality, we note that
%
%
\begin{equation}\label{1.27}
\Phi_{y,L}(0) = \varphi_{y,L} - \varphi_{0,L} + P_y[H_0 > T_{U_L}]
\varphi_{0,L}.
\end{equation}
By the strong Markov property and the symmetry of $g_L(\cdot,\cdot)$
one has
\[
P_y[H_0 > T_{U_L}] = \bigl(g_L(0,0) - g_L(0,y)\bigr)
/ g_L(0,0) \stackrel{\mbox{\fontsize{8.36}{8.36}\selectfont{(\ref
{1.26})}}}{\le} c(y) / g_L(0,0),
\]
and by~(\ref{1.20}) one finds that
%
%
\begin{equation}\label{1.28}
E[\varphi^2_{0,L}] = g_L(0,0) .
\end{equation}
Since $\lim_L g_L(0,0) = \infty$, it follows that the last term of
(\ref{1.27}) converges to $0$ in $L^2$ as $L$ tends to infinity, and
the second equality of~(\ref{1.23}) now follows.
\end{pf}

We thus introduce on some auxiliary probability space
%
%
\begin{equation}\label{1.29}
\begin{tabular}{p{305pt}}
$\psi_y, y \in\IZ^2$, a centered Gaussian field with covariance
function
$E[\psi_y \psi_{y^\prime}] = 3(a(y) + a(y^\prime)
- a(y^\prime-y)), y,y^\prime\in\IZ^2$.
\end{tabular}\hspace*{-32pt}
\end{equation}
Up to an inessential multiplicative factor $\sqrt{3}$, we can thus
interpret $\psi_y, y \in\IZ^2$, as the field of ``increments at the
origin'' of the two-dimensional massless free field, or as the
two-dimensional massless free field pinned at the origin.

The last topic of this section concerns zero-dimensional Bessel
processes. We denote by $\mathrm{B E S}^0(a,\tau)$ the law at time $\tau\ge0$
of a zero-dimensional Bessel process starting from $a \ge0$. If $R$ is
a random variable with distribution $\mathrm{BES}^0(a,\tau)$, the Laplace
transform of $R^2$ is given by the formula (see~\cite{RevuYor98},
page 411, or~\cite{IkedWata89}, pa\-ge~239)
%
%
\begin{equation}\label{1.30}
E[e^{-\lambda R^2}] = \exp\biggl\{ - \frac{\lambda a^2}{1 + 2
\tau\lambda}\biggr\}\qquad \mbox{for } \lambda\ge0 .
\end{equation}
We also denote by $\operatorname{BESQ}^0(a^2,\tau)$ the law of $R^2$; this is the
distribution of a zero-dimensional square Bessel process at time
$\tau$, starting from $a^2$ at time~$0$.

\section{Laplace functional of occupation times}\label{sec2}

In this section we obtain a formula for the Laplace functional of the
occupation times $L_{x,u}$, which proves~(\ref{0.3}); see Theorem \ref
{theo2.3}. As a by-product we note the absence for fixed $u$ of a
global factorization for the field $L_{x,u} - L_{0,u}$, $x \in\IZ^d$,
similar to that of the limit law in~(\ref{0.8}), even through each
individual variable $L_{x,u} - L_{0,u}$ is distributed as the product
of a time-marginal of a zero-dimensional Bessel process with an
independent centered Gaussian variable; see Remark~\ref{rem2.4}(2).
Preparatory Theorem~\ref{theo2.1} will be repeatedly used in the
sequel and shows in particular that the characteristic function of a
finite linear combination of the variables $L_{x,u}$, $x \in\IZ^d$,
is analytic in the neighborhood of the origin. This will play an
important role in Section~\ref{sec4}.

We denote by $G$ the linear operator
%
%
\begin{equation}\label{2.1}
Gf(x) = \sum_{x^\prime\in\IZ^d} g(x,x^\prime) f(x'),\qquad x \in
\IZ^d ,
\end{equation}
which is well defined when $\sum_{x^\prime} g(x,x')|f(x')| < \infty$,
and in particular when $f$ vanishes outside a finite set. When $V$ is a
function on $\IZ^d$ vanishing outside a finite set, we write $GV$ for
the composition of $G$ with the multiplication operator by~$V$, so that
$GV$ naturally operates on $L^\infty(\IZ^d)$ (we recall that \mbox{$\|
\cdot
\|_\infty$} denotes the corresponding sup-norm; see the beginning of
Section~\ref{sec1}).

\begin{theorem}\label{theo2.1}
If $V$ has support in $K \subset\subset\IZ^d$, then for $\|V\|
_\infty\le c(K)$,
%
%
\begin{equation}\label{2.2}
\|GV \|_{L^\infty\rightarrow L^\infty} < 1
\end{equation}
and for any $u \ge0$,
%
%
\begin{equation}\label{2.3}
\IE\biggl[\exp\biggl\{\sum_{x \in\IZ^d} V(x)L_{x,u}\biggr\}\biggr] =
\exp\bigl\{ u\bigl(V, (I-GV)^{-1} 1\bigr)\bigr\}.
\end{equation}
\end{theorem}
\begin{pf}
Claim\vspace*{1pt}~(\ref{2.2}) is immediate. As for~(\ref{2.3}), we note that
defining for $x \in\IZ^d$, $u \ge0$, the function on $\widehat
{W}^*$ [see above~(\ref{0.1}) for notation]
\[
\gamma_x(\widehat{w}^*) = \sum_{n \in\IZ} \sigma_n(\widehat{w})
1\{X_n (\widehat{w}) = x\}\qquad \mbox{for any $\widehat{w} \in
\widehat{W}$}\qquad\mbox{with $\pi^*(\widehat{w}) = \widehat{w}^*$},
\]
we have the identity
%
%
\begin{eqnarray}\label{2.4}
L_{x,u}(\omega) = \sum_i \gamma_x(\widehat{w}^*_i) 1\{u_i \le
u\}\nonumber\\[-8pt]\\[-8pt]
&&\eqntext{\mbox{for } \omega=\displaystyle  \sum_i \delta_{(\widehat{w}^*_i,u_i)} \in
\Omega, x \in\IZ^d, u \ge0 .}
\end{eqnarray}
The interlacement point process $\omega$ is Poisson with intensity measure
$\widehat{\nu}(d\widehat{w}^*) \,du$ under $\IP$, and hence when $V$ is
supported in $K \subset\subset\IZ^d$ and $\|V\|_\infty\le c(K)$, we
have
%
%
\begin{eqnarray}\label{2.5}
&&\IE\biggl[\exp\biggl\{\sum_{x \in\IZ^d} V(x)L_{x,u}\biggr\}\biggr]
\nonumber\\
&&\qquad= \IE\biggl[\exp\biggl\{\int_{\widehat{W}^* \times\IR_+}
\sum_{x \in\IZ^d} V(x) \gamma_x(\widehat{w}^*) 1\{v \le u\}\,
d\omega(\widehat{w}^*,v)\biggr\}\biggr]
\nonumber\\
&&\qquad= \exp\biggl\{ u \int_{\widehat{W}^*} \bigl(e ^{\sum
_{x \in\IZ^d} V(x) \gamma_x(\widehat{w}^*)} - 1\bigr) \,d \widehat
{\nu}(\widehat{w}^*)\biggr\}
\\
&&\qquad\stackrel{\mbox{\fontsize{8.36}{8.36}
\selectfont{(\ref{0.1})}}}{=} \exp\bigl\{ u E_{e_K} \bigl[e ^{\sum
_{x \in\IZ^d} V(x) \sum_{k \ge0} \sigma_k 1\{X_k
= x\}}-1\bigr]\bigr\}
\nonumber\\
&&\qquad= \exp\biggl\{ u E_{e_K} \biggl[ \exp\biggl\{\int^\infty_0
V(\overline{X}_s) \,ds\biggr\} - 1\biggr]\biggr\},\nonumber
\end{eqnarray}
where for $s \ge0$, $\widehat{w} \in\widehat{W}_+$, we have set
\[
\overline{X}_s(\widehat{w}) = X_k(\widehat{w})\qquad \mbox{when }
\sigma_0(\widehat{w}) +\cdots+ \sigma_{k-1}(\widehat{w}) \le s <
\sigma_0(\widehat{w}) +\cdots+ \sigma_k(\widehat{w})
\]
(by convention the term bounding $s$ from below vanishes when $k=0$),
that is, $\overline{X}_\point$ is the natural continuous time random
walk on $\IZ^d$ with jump parameter $1$ defined on $\widehat{W}_+$.
Thus for $\|V\|_\infty\le c(K)$ we find by a classical calculation that
\begin{eqnarray}\label{2.6}
&&
E_{e_K}\bigl[e^{\int_0^\infty V(\overline{X}_s)\,ds}\bigr] \nonumber\\
&&\qquad= E_{e_K}
\biggl[\sum_{n \ge0} \frac{1}{n!} \biggl(\int_0^\infty
V(\overline{X}_s) \,ds\biggr)^n\biggr]
\nonumber\\[-8pt]\\[-8pt]
&&\qquad=\sum_{n \ge0} E_{e_K}\biggl[\int_{0 < s_1 <\cdots< s_n < \infty
} V(\overline{X}_{s_1}) \cdots V(\overline{X}_{s_n}) \,ds_1 \cdots
ds_n\biggr]
\nonumber\\
&&\qquad=\operatorname{cap} (K) + \sum_{n \ge1} \sum_{x \in K} e_K(x) [(GV)^n
1](x),\nonumber
\end{eqnarray}
using Fubini's theorem and the Markov property in the last step. Since
$V$ vanishes outside $K$, it also follows from~(\ref{1.5}) that for $n
\ge1$, one has
%
%
\begin{equation}\label{2.7}
\sum_x e_K(x) [(GV)^n 1] (x) = \sum_{x'} V(x')[(GV)^{n-1}1](x') =
(V,(GV)^{n-1} 1).\hspace*{-28pt}
\end{equation}
As a result we see that when $\|V\|_\infty\le c(K)$,
%
%
\begin{equation}\label{2.8}
E_{e_K}\bigl[e^{\int_0^\infty V(\overline{X}_s)\,ds} -1\bigr] = \sum
_{n \ge1} (V,(GV)^{n-1}1) = \bigl(V,(I-GV)^{-1}1\bigr).
\end{equation}
Inserting this identity in the last line of~(\ref{2.5}) completes the
proof of Theorem~\ref{theo2.1}.
\end{pf}
\begin{remark}\label{rem2.2}
As a staightforward consequence of Theorem~\ref{theo2.1}, see that for
any finitely supported real valued function $V$ on $\IZ^d$ and $u \ge
0$, the random variable $\sum_{x \in\IZ^d} V(x) L_{x,u}$ has a
characteristic function which coincides in the neighborhood of the
origin with the exponential of an analytic function. So this
characteristic function is analytic in the sense of Chapter 7 of~\cite{Luka70}.
Particularly with Theorem 7.1.1, page 193\vadjust{\goodbreak} of~\cite{Luka70},
one has the identity
%
%
\begin{equation}\label{2.9}
\IE\biggl[\exp\biggl\{ z \sum_x V(x) L_{x,u}\biggr\}\biggr] = \Phi
_{V,u}(z),\qquad z \in S,
\end{equation}
where in the above formula $S$ stands for the maximal vertical strip in
$\IC$ to which the function $z \rightarrow\exp\{u \sum_{n \ge1}$
$z^n(V,(GV)^{n-1} 1)\}$ can be\vspace*{1pt} an analytically extended, and $\Phi
_{V,u}$ for this extension [i.e., for $z \in S$, $\exp\{z \sum_x
V(x) L_{x,u}\}$ is integrable, and the equality~(\ref{2.9}) holds].
This fact will be very helpful and repeatedly used in the sequel.
\end{remark}

We now derive an alternative expression for the right-hand side of
(\ref{2.3}) and need some additional notation for this purpose.

For $I \subset\subset\IZ^d$ nonempty, we denote by $G_I$ the matrix
$g(x,x')$, $x,x' \in I$. It is well known to be positive definite (see,
e.g., Lemma 3.3.6 of~\cite{MarcRose06}), and we introduce
%
%
\begin{equation}\label{2.10}
g_I = \operatorname{det} (G_I) > 0,
\end{equation}
where the right-hand side does not depend on the identification of $I$
with $\{1,\ldots,|I|\}$ we use. We also set by convention $g_I = 1$,
when $I = \phi$. Further we introduce
%
%
\begin{equation}\label{2.11}\qquad
\mbox{$c_I = \mbox{the}$ sum of all coefficients of the matrix of cofactors
of $G_I$},
\end{equation}
and note that $c_I$ does not depend on the identification of $I$ with
$\{1,\ldots, | I |\}$ we employ; for instance, $c_I/g_I$ coincides with
the sum of all coefficients of the inverse matrix of $G_I$. The above
also shows that
%
%
\begin{equation}\label{2.12}
c_I > 0\qquad\mbox{when $I \subset\subset\IZ^d$ is nonempty}.
\end{equation}
We extend the notation to the case $I = \phi$ with the convention
$c_\phi= 0$. It is known (see~\cite{Spit01}, page 301) that
%
%
\begin{equation}\label{2.13}
\operatorname{cap}(I) = c_I / g_I\qquad \mbox{for all $I \subset\subset\IZ^d$}.
\end{equation}
We are now ready for the main result of this section.
\begin{theorem}\label{theo2.3}
When $V$ has support in $K \subset\subset\IZ^d$, then for $\|V\|
_\infty\le c(K)$, $u \ge0$,
%
%
\begin{equation}\label{2.14}
\IE\biggl[\exp\biggl\{ - \sum_x V(x) L_{x,u}\biggr\}\biggr] = \exp
\biggl\{ - u \frac{\sum_{I \subset K} c_I V_I}{\sum_{I
\subset K} g_I V_I} \biggr\},
\end{equation}
where $V_I \stackrel{\mathrm{def}}{=} \Pi_{x \in I} V(x)$ ($V_I =1$, by
convention when $I = \phi$).

In addition~(\ref{2.14}) holds whenever $V$ is nonnegative and
vanishes outside~$K$.
\end{theorem}
\begin{pf}
By Theorem~\ref{theo2.1} we know that for $\|V\|_\infty\le c(K)$, the
left-hand side of~(\ref{2.14}) equals $\exp\{- u (V, (I+GV)^{-1} 1)\}
$. With no loss of generality we assume that $|K| \ge2$. We
identify\vadjust{\goodbreak}
$K$ with $\{1,\ldots,n\}$, where $n = |K|$, via an enumeration
$x_1,\ldots,x_n$ of $K$. Writing $v_\ell= V(x_\ell)$, we see that
%
%
\begin{equation}\label{2.15}
\bigl(V,(I + GV)^{-1} 1\bigr) = \sum^n_{k,\ell= 1} \pmb{C}_{k,\ell} v_\ell
\big/ \operatorname{det} (\mathbf{I} + \pmb{GV}) ,
\end{equation}
where $\pmb{G}$ stands for the $n \times n$ matrix $g(x_k,x_\ell)$,
$1 \le k, \ell\le n$, $\pmb{V}$ for the diagonal matrix with
coefficients $v_\ell$, $1 \le\ell\le n$, on the diagonal, $\pmb{I}$
for the identity matrix and $\pmb{C}$ for the matrix of cofactors of
$\pmb{I} + \pmb{GV}$. Observe that for $1 \le k, \ell\le n$, one has
%
%
\begin{equation}\label{2.16}
\pmb{C}_{k,\ell} = \operatorname{det}\bigl((\pmb{I} + \pmb{GV})^{k,\ell}\bigr),
\end{equation}
where $(\pmb{I} + \pmb{GV})^{k,\ell}$ stands for the $n \times n$
matrix, where the $k$th line and the $\ell$th column of $\pmb{I} +
\pmb{GV}$ have been replaced by $0$, except for the coefficient at
their intersection, which is replaced by $1$.

Given an $n \times n$ matrix $\pmb{A} = (a_{k,\ell})$, we develop the
determinant of $\pmb{A}$ according to the classical formula
%
%
\begin{equation}\label{2.17}
\operatorname{det} \pmb{A} = \sum_\sigma\operatorname{sign}(\sigma) \prod
^n_{k=1} a_{k,\sigma(k)},
\end{equation}
where $\sigma$ runs over the permutations of $\{1,\ldots,n\}$, and
$\operatorname{sign}(\sigma)$ denotes the signature of $\sigma$.

We now develop the determinant $\operatorname{det} (\pmb{I} + \pmb{GV})$. For
each subset $J \subseteq\{1,\ldots,n\}$ we collect the terms
corresponding to permutations $\sigma$ of $\{1,\ldots,n\}$ such that
$\sigma(k) = k$, for $k \in J$, the choice of $1$ in each term
$(1+g(0) v_k)$, $k \in J$, and for any $k \notin J$ such that $\sigma
(k) = k$, the choice of $g(0) v_k$ instead. Thus for each such~$J$,
setting $\wt{J} = \{1,\ldots,n\} \setminus J$, the sum of these terms
equals $\operatorname{det} (\pmb{G}_{| \wt{J} \times\wt{J}}) \prod_{\ell
\in\wt{J}} v_\ell$. Thus summing over all subsets $J$ of $\{1,\ldots
,n\}$ we find
%
%
\begin{equation}\label{2.18}
\operatorname{det}(\pmb{I} + \pmb{GV}) = \sum_{I \subseteq K} g_I V_I .
\end{equation}
We now turn to the numerator of the right-hand side of~(\ref{2.15}).
We use the convention $\{k,\ell\} = \{k\} = \{\ell\}$, when $k = \ell
$. As above we develop the determinant $\operatorname{det}((\pmb{I} + \pmb
{GV})^{k,\ell})$; see~(\ref{2.16}),~(\ref{2.17}). We can assume that
the permutations $\sigma$ of $\{1,\ldots,n\}$ entering the development
satisfy $\sigma(k) = \ell$. For each $J \subseteq\{1,\ldots,n\}
\setminus\{k,\ell\}$, we collect the terms corresponding to
permutations $\sigma$ such that $\sigma(m) = m$, for $m \in J$, the
choice of $1$ in each term $(1 + g(0) v_m)$, for $m \in J$, and for any
$m \notin J \cup\{k,\ell\}$ with $\sigma(m) = m$, the choice of
$g(0) v_m$ instead. Setting $\wt{J} = \{1,\ldots,n\} \setminus J$, we
see that the sum of these terms for a fixed given $J$ as above equals
$\Pi_{m \in\wt{J} \setminus\{\ell\}} v_m \operatorname{det}(\pmb
{G}_{\wt{J} \times\wt{J}}^{k,\ell})$, where $\pmb{G}_{\wt{J}
\times\wt{J}}^{k,\ell}$ stands for the matrix where the $k$th line
and the $\ell$th column of the matrix $\pmb{G}_{| \wt{J} \times\wt
{J}}$ (i.e., $\pmb{G}$ restricted to $\wt{J} \times\wt{J}$) are
replaced by zero except for the coefficient at their intersection,
which is replaced by $1$. Thus summing over all possible $J \subseteq\{
1,\ldots, n\} \setminus\{k,\ell\}$ and all $k, \ell$ in $\{1,\ldots
,n\}$, we obtain
%
%
\begin{eqnarray}\label{2.19}
\sum_{k,\ell= 1}^n \pmb{C}_{k,\ell} v_\ell& = & \sum^n_{k,\ell=
1} \sum_{H \supseteq\{k,\ell\}} \prod_{m \in H} v_m
\operatorname{det} (\pmb{G}^{k,\ell}_{H \times H})
\nonumber\\
& = &\sum_{\phi\not= H \subseteq\{1,\ldots,n\}} \prod_{m
\in H} v_m \sum_{k,\ell\in H} \operatorname{det} (\pmb{G}^{k,\ell
}_{H \times H})
\\
&\stackrel{\mbox{\fontsize{8.36}{8.36}\selectfont{(\ref{2.11})}}}{=}&
\sum_{I \subseteq K} V_I
c_I
\nonumber
\end{eqnarray}
(using the convention $c_\phi= 0$ in the last equality).

Combining~(\ref{2.15}),~(\ref{2.16}),~(\ref{2.18}) we obtain (\ref
{2.14}). Finally in the case of a nonnegative $V$ with support in $K$,
we note that $\IE[\exp\{ - z \sum_x V(x) L_{x,u}\}]$ is analytic in
the strip $\operatorname{Re} z > 0$, and coincides for small positive $z$ with the
function $\exp\{ - u \sum_{I \subseteq K} c_I V_I z^{|I|} /\sum
_{I \subseteq K} g_I V_I z^{|I|}\}$, which is analytic in the
neighborhood of the positive half-line. Both functions thus coincide
for $z=1$, and our last claim follows.
\end{pf}
\begin{remark}\label{rem2.4}
(1) Choosing $V = \lambda1_K$, with $\lambda\ge0$ and $K \subset
\subset\IZ^d$, we deduce from~(\ref{2.14}) by letting $\lambda$
tend to infinity that
\[
\IP[L_{x,u} = 0\mbox{, for all } x \in K] = \exp\biggl\{- u
\frac{c_K}{g_K}\biggr\}\qquad \mbox{for $u \ge0$}.
\]
Introducing the interlacement at level $u$,
%
\begin{eqnarray}
\cI^u (\omega) = \{x \in\IZ^d\mbox{; for some $i \ge0$ such that
$u_i \le u$, $\widehat{w}^*_i$ enters $x$}\}\nonumber\\
&&\eqntext{\mbox{if }\displaystyle \omega= \sum_{i \ge0}
\delta_{(\widehat{w}^*_i,u_i)},}
\end{eqnarray}
and taking~(\ref{2.13}) into account, we recover the well-known
formula (see~(\ref{2.16}) of~\cite{Szni10a})
%
%
\begin{equation}\label{2.20}
\IP[\cI^u \cap K = \phi] = \exp\{ - u \operatorname{cap} (K)\}\qquad
\mbox{for } u \ge0, K \subset\subset\IZ^d .
\end{equation}

(2) Choosing $V= \lambda1_{\{x\}}$, with $\lambda\ge0$ and $x \in
\IZ^d$, Theorem~\ref{theo2.3} now yields that
%
%
\begin{equation}\label{2.21}
E [\exp\{ - \lambda L_{x,u}\} ] = \exp\biggl\{ - \frac{\lambda
u}{1 + g(0) \lambda}\biggr\}\qquad \mbox{for $\lambda\ge0$},
\end{equation}
and in view of~(\ref{1.30}) we find that
%
%
\begin{equation}\label{2.22}
L_{x,u}\mbox{ is }\operatorname{BESQ}^0\biggl(u, \frac{g(0)}{2}\biggr)\mbox{
distributed.}
\end{equation}
If $x,x' \in\IZ^d$ are distinct, choosing $V = z(1_{\{x\}} - 1_{\{x'\}
})$ in~(\ref{2.14}) with $z$ small and real and extending the identity
to $z = it$, $t \in\IR$, with the help of~(\ref{2.9}), we find that
\begin{eqnarray}\label{2.23}
&&
E[\exp\{ it(L_{x',u} - L_{x,u})\}] \nonumber\\[-8pt]\\[-8pt]
&&\qquad= \exp\biggl\{ -2 u \frac
{(g(0) - g(x' -x))t^2}{1 + (g(0)^2 - g(x' -x)^2)t^2}\biggr\}\qquad
\mbox{for $t \in\IR$}.\nonumber
\end{eqnarray}
In view of~(\ref{1.30}) we thus find that
%
%
\begin{equation}\label{2.24}
\begin{tabular}{p{315pt}}
$L_{x,u} - L_{x',u}$ has the law of $R \psi$, where $R$ and $\psi
$ are independent, respectively, $\mathrm{BES}^0 ( \sqrt{u}, \frac{g(0) +
g(x-x')}{4})$, and centered
Gaussian with variance $4(g(0) - g(x-x'))$ distributed.
\end{tabular}\hspace*{-32pt}
\end{equation}
Let us, however, point out that in the case of three distinct points
$x,x',x''$ in $\IZ^d$, the law of the random vector $(L_{x',u} -
L_{x,u}, L_{x'',u} - L_{x,u})$ does not coincide with that of the
scalar multiplication of a two-dimensional Gaussian vector by an
independent $\mathrm{BES}^0(a,\tau)$-variable, when $u > 0$. Indeed one has
\[
\IP[L_{x',u} - L_{x,u} > 0, L_{x'',u} - L_{x,u} = 0] \ge\IP[\cI^u
\ni x', \cI^u \cap\{x,x''\} = \phi] > 0
\]
as a consequence of~(\ref{2.20}) and the fact that $\operatorname{cap}(\{
x,x',x''\}) > \operatorname{cap}(\{x,x''\})$. But for the above mentioned
distribution both components necessarily vanish simultaneously on a set
of full measure, and the above probability would equal zero if such an
identity in law was to hold. We will, however, see in Section~\ref{sec5} how
$L_{x,u}$, $x \in\IZ^d$, can be related to the $d$-dimensional
Gaussian free field, by letting $u$ tend to infinity, instead of
keeping $u$ fixed.

(3) Random interlacements can be related to the Poissonian gas of
Markov loops; see~\cite{Leja10,Leja11}. Heuristically they
correspond to ``loops passing through infinity;'' see~\cite{Leja11},
page 85. The identity for Markov loops corresponding to~(\ref{2.3}) of
Theorem~\ref{theo2.1} above can be found in Corollary 1 of Chapter 4,
Section 1 and Proposition 7 of Chapter 2, Section 4 of~\cite{Leja11}.
The presence of a logarithm and a trace in the expressions leading to
Proposition 7 of~\cite{Leja11} is emblematic of the Markov loop
measure and can be contrasted with the expression in~(\ref{2.8}) for
random interlacements [which is then inserted in the last line of (\ref
{2.5})]. In the case of a Poissonian gas of Markov loops on a finite
weighted graph with a suitable killing, it is shown in Theorem 13 of
\cite{Leja10} that the occupation field of the gas of loops at level
$\frac{1}{2}$ (playing the role of $u$ in the context of
\cite{Leja10}) is distributed as half the square of a centered Gaussian
free field with covariance the corresponding Green density. For similar
reasons as in (2) above, no such identity holds for random
interlacements at any fixed level $u$. We will, however, present in the
next two sections limiting procedures that relate random interlacements
to the Gaussian free field.

(4) As we now explain, the results of this section can be extended to
the case of continuous time random interlacements on a transient
weighted graph. One considers a countable connected graph $E$ which is
locally finite and endowed with nonnegative symmetric weights $\rho
_{x,x'} = \rho_{x',x}$, which are positive exactly when $\{x,x'\}$
belongs to the edge set $\cE$ of $E$. One assumes that the induced
random walk with transition probability $p_{x,x'} = \rho_{x,x'} / \rho
(x)$, where $\rho(x) = \sum_{x' \in E} \rho_{x,x'}$, is transient.
Random interlacements can be constructed on such a transient weighted
graph; see~\cite{Szni10a}, Remark 1.4, and~\cite{Teix09b}. Continuous
time random interlacements can also be constructed, in essence by the
same procedure described in the \hyperref[intro]{Introduction}, endowing the discrete
doubly infinite paths with i.i.d. exponential variables of parameter
$1$. The corresponding expression for the measure $\wh{Q}_K$, for $K$
finite subset of $E$, remains the same as in~(\ref{0.1}); simply, the
expression for $e_K(\cdot)$ the equilibrium measure of $K$, which
appears in~(\ref{1.4}), now has to be multiplied by the factor $\rho
(x)$ in the present context.

The occupation time variables $L_{x,u}$, $x \in E$, $u \ge0$, are
defined by a similar formula as in~(\ref{0.2}), but the expression on
the right-hand side of~(\ref{0.2}) is now divided by $\rho(x)$. The
linear operator $G$ corresponding to~(\ref{2.1}) operates, say, on
functions $f$ on $E$ with finite support, via the formula
\[
Gf(x) = \sum_{x' \in E} g(x,x') f(x') \rho(x'),\qquad x \in E,
\]
where $g(\cdot,\cdot)$ now stands for the Green density, which is
obtained by dividing the expression corresponding to the right-hand
side of~(\ref{1.1}) by $\rho(x')$.

The proof of Theorem~\ref{theo2.1} can be adapted to this context to
show that when $K$ is a finite subset of $E$, and $V$ has support in
$K$, then for $\|V\|_{L^\infty(E)}$ sufficiently small, $\|GV\|
_{L^\infty(E) \rightarrow L^\infty(E)} < 1$, and for any $u \ge0$,
%
%
\begin{equation}\label{2.25}
\IE\biggl[\exp\biggl\{ \sum_{x \in E} V(x) L_{x,u} \rho(x)\biggr\}
\biggr] = \exp\bigl\{ u\bigl(V, (I-GV)^{-1} 1\bigr)\bigr\},
\end{equation}
where now $(f,g)$ stands for $\sum_{x \in E} f(x) g(x) \rho(x)$
(whenever this sum is absolutely convergent). Likewise the proof of
Theorem~\ref{theo2.3} is easily adapted, and one finds that for $V$ as
above, $u \ge0$,
%
%
\begin{equation}\label{2.26}
\IE\biggl[\exp\biggl\{ -\sum_{x \in E} V(x) L_{x,u} \biggr\}\biggr] =
\exp\biggl\{ - u \frac{\sum_{I \subseteq K} c_I \Pi_{x \in
I} V(x)}{\sum_{I \subseteq K} g_I \Pi_{x \in I} V(x)}\biggr\}
\end{equation}
with $g_I$ and $c_I$ defined as in~(\ref{2.10}),~(\ref{2.11}) [with
$g(\cdot,\cdot)$ now denoting the Green density].

(5) One can define the stationary field of discrete occupation times
$\ell_{x,u}$, $x \in\IZ^d$, $u \ge0$, analogously to $L_{x,u}$,
simply replacing\vspace*{1pt} $\sigma_n$ by $1$ in~(\ref{0.2}). When $V$ is a
function on $\IZ^d$ with support contained in $K \subset\subset\IZ
^d$, it follows that $1 - e^{-V}$ is a function supported in $K$ with
values in $(-\infty,1)$, and one has the identity
%
%
\begin{equation}\label{2.27}
\IE\biggl[\exp\biggl\{\sum_x V(x) \ell_{x,u}\biggr\}\biggr] = \IE
\biggl[\exp\biggl\{\sum_x \bigl(1 - e^{-V(x)}\bigr) L_{x,u}\biggr\}\biggr],
\qquad u \ge0,\hspace*{-35pt}
\end{equation}
as can be seen by integrating out the exponential variables in the
right member of~(\ref{2.27}) (of course both members of the above
equality may be infinite). As a result, Theorems~\ref{theo2.1} and
\ref{theo2.3} also yield identities concerning the Laplace
functional of $(\ell_{x,u})_{x \in\IZ^d}$.
\end{remark}

\section{Preparation for the study of long rods}\label{sec3}

In this section we introduce notation specific to $\IZ^3$ and provide
estimates in Lemmas~\ref{lem3.1} and~\ref{lem3.2}, which will be
recurrently used in the next section, when we investigate the
occupation times spent by interlacements at a suitably scaled level in
long rods. These controls will play an important role in the asymptotic
analysis of the power series entering the characteristic functions of
these occupation times. Throughout this section we assume that $d=3$,
and constants depend on the finite subset $\Lambda$ of $\IZ^2$
introduced in~(\ref{3.1}) below. The notation $\| \cdot\|_\infty$
refers to the supremum norm $\| \cdot\|_{L^\infty(B)}$, where $B$
appears in~(\ref{3.1}).

We consider $\Lambda\subset\subset\IZ^2$ containing $0$ and $N >
1$. We also define
%
%
\begin{equation}\label{3.1}
B = \Lambda\times J\qquad \mbox{where } J = \{1, \ldots, N\} .
\end{equation}
We write $\pi_{\IZ^2}$ and $\pi_\IZ$ for the respective $\IZ^2$-
and $\IZ$-projections on $\IZ^3$ identified with $\IZ^2 \times\IZ
$. Given a function $F$ on $B$, we write $\langle F \rangle$ for the
function obtained by averaging $F$ on horizontal layers and $\langle F
\rangle_z$ for the average of $F$ on the layer $\Lambda\times\{z\}$,
so that
%
%
\begin{equation}\label{3.2}\quad
\langle F \rangle(x) = \langle F \rangle_z = \frac{1}{|\Lambda
|} \sum_{y \in\Lambda} F((y,z))\qquad \mbox{for $x \in B$ with $\pi
_\IZ(x) = z$}.
\end{equation}
We also introduce the function
%
%
\begin{equation}\label{3.3}
[F]_0(x) = F((0,z))\qquad\mbox{for $x$ in $B$ with $\pi_\IZ(x) = z$}.
\end{equation}
It is plain that for any function $F$ on $B$,
%
%
\begin{equation}\label{3.4}
\bigl\langle F - \langle F \rangle\bigr\rangle= 0
\end{equation}
and that
%
%
\begin{equation}\label{3.5}
F= \langle F \rangle\qquad\mbox{when $F$ only depends on the $\IZ$-component.}
\end{equation}
We consider nonempty sub-intervals of $J$,
%
%
\begin{equation}\label{3.6}
I_0 \subseteq I_1 \subsetneq J\qquad \mbox{with } L = d(I_0, J
\setminus I_1) \ge1
\end{equation}
(we refer to the beginning of Section~\ref{sec1} for notation). We write
%
%
\begin{equation}\label{3.7}
C_0 = \Lambda\times I_0 \subseteq C_1 = \Lambda\times I_1 \subseteq B.
\end{equation}
We recall the convention concerning constants and the notation \mbox{$\|
\cdot\|_\infty$} stated at the beginning of this section. The
estimates in the next lemma reflect the decay at infinity of the Green
function [see~(\ref{1.12})] and the fact that the discrete gradient of
$g(\cdot)$ has an improved decay at infinity; see~(\ref{3.15}) below.
\begin{lemma}\label{lem3.1}
For any function $F$ on $B$, one has:
%
%
\begin{eqnarray}
\label{3.8}
\|GF\|_\infty&\le& c\log N \|F\|_\infty;
\\
\label{3.9}
\|1_{C_0} G 1_{B \setminus C_1} F\|_\infty&\le& c \log\biggl(
\frac{N+1}{L}\biggr) \|F\|_\infty;
\\
\label{3.10}
\|GF\|_\infty&\le& c \|F\|_\infty\qquad \mbox{when } \langle F
\rangle=0;
\\
\label{3.11}
\|1_{C_0} G 1_{B \setminus C_1} F\|_\infty&\le&\frac{c}{L}
\|1_{B \setminus C_1} F\|_\infty\qquad \mbox{when } \langle F \rangle
= 0;
\\
\label{3.12}
\|GF - [GF]_0\|_\infty&\le& c \|F\|_\infty;
\\
\label{3.13}
\hspace*{36pt}\| 1_{C_0} (G1_{B\setminus C_1} F - [G1_{B\setminus C_1}F]_0 )\|
_\infty&\le&\frac{c}{L} \|1_{B\setminus C_1} F\|_\infty.
\end{eqnarray}
\end{lemma}
\begin{pf}
We begin with~(\ref{3.9}) and note that for $x \in C_0$,
\begin{eqnarray*}
|(G1_{B\setminus C_1}F)(x) | &=& \biggl| \sum_{x' \in B\setminus C_1}
g(x,x') F(x')\biggr|
\stackrel{\mbox{\fontsize{8.36}{8.36}\selectfont{(\ref{1.3})}}}{\le}
c\biggl(\sum_{L \le k \le N} \frac{1}{k}\biggr) \|F\|_\infty\\
&\le& c \log
\biggl(\frac{N+1}{L}\biggr) \|F\|_\infty,
\end{eqnarray*}
whence~(\ref{3.9}). The bound~(\ref{3.8}) is proved in the same fashion.

We then turn to the proof of~(\ref{3.11}) and note that when $\langle
F\rangle= 0$, for $x \in C_0$ one has with the notation $x' =
(y',z')$, $\overline{x} = (\overline{y},z')$ \mbox{[so that $\pi_\IZ(x')
= \pi_\IZ(\overline{x}) = z'$]}
%
%
\begin{eqnarray}\label{3.14}
(G 1_{B\setminus C_1}F) (x) & = &\sum_{x' \in B \setminus C_1}
g(x,x') F(x') = \sum_{z' \in J \setminus I_1} \sum_{y' \in
\Lambda} g(x,x') F(x')
\nonumber\\
& \stackrel{\langle F \rangle= 0}{=}& \frac{1}{|\Lambda
|} \sum_{z' \in J \setminus I_1} \sum_{y', \overline{y} \in
\Lambda} g(x,x')\bigl(F(x') - F(\overline{x})\bigr)
\\
& = &\frac{1}{|\Lambda|} \sum_{z' \in J \setminus I_1}
\sum_{y', \overline{y} \in\Lambda} \bigl(g(x,x') - g(x,\overline{x})\bigr)
F(x') .
\nonumber
\end{eqnarray}
From Theorem 1.5.5, page 32 of~\cite{Lawl91}, one knows that
%
%
\begin{equation}\label{3.15}
|g(x+a) - g(x)| \le\frac{c |a|}{1 + |x|^2} \qquad\mbox{for
$x \in\IZ^3$, $|a| \le\operatorname{diam} (\Lambda)$},
\end{equation}
where $\operatorname{diam} (\Lambda)$ stands for the diameter of $\Lambda$. As
a result we see that
%
%
\begin{equation}\label{3.16}
|(G1_{B \setminus C_1}F)(x)| \le c \|1_{B \setminus C_1} F\|_\infty
\sum_{k \ge L} \frac{1}{1 + k^2} \le\frac{c}{L} \|
1_{B \setminus C_1} F\|_\infty,
\end{equation}
whence~(\ref{3.11}). One proves~(\ref{3.10}) analogously.

As for~(\ref{3.13}) we note that with the notation $x_0 = (0,z)$ when
$x = (y,z)$, we have for $x \in C_0$,
%
%
\begin{eqnarray}\label{3.17}
|(G1_{B \setminus C_1} F- [G1_{B \setminus C_1} F]_0)(x)| & = &\biggl|
\sum_{x'\in B \setminus C_1} \bigl(g(x,x') - g(x_0,x')\bigr) F(x')\biggr|
\nonumber\\
& \stackrel{\mbox{\fontsize{8.36}{8.36}\selectfont{(\ref{3.15})}}}{\le}&
c \sum_{L \ge k}
\frac{1}{1 + k^2} \|1_{B \setminus C_1} F\|_\infty\\
&\le&\frac
{c}{L} \|1_{B \setminus C_1} F\|_\infty,
\nonumber
\end{eqnarray}
whence~(\ref{3.13}). The bound~(\ref{3.12}) is proved analogously.
\end{pf}

We conclude this section with the following lemma.
\begin{lemma}\label{lem3.2}
For $F,H$ functions on $B$, one has
%
%
\begin{equation}\label{3.18}
\| \langle F(GH)\rangle\|_\infty\le c \|F\|_\infty\| H\|_\infty
\qquad\mbox{when } \langle F\rangle= 0 \mbox{ or } \langle H
\rangle= 0.
\end{equation}
\end{lemma}
\begin{pf}
The case $\langle H \rangle= 0$ is immediate thanks to~(\ref{3.10}).
In the case where $\langle F \rangle= 0$, we write for $z \in J$, with
the notation $x = (y,z)$, $\overline{x} = (\overline{y},z)$,
%
%
\begin{eqnarray}\label{3.19}
\langle F(GH)\rangle_z & = &\frac{1}{|\Lambda|} \sum_{y \in
\Lambda} F(x) \sum_{x'\in B} g(x,x') H(x')
\nonumber\\
& \stackrel{\langle F\rangle= 0}{=} &\frac{1}{|\Lambda
|^2} \sum_{y, \overline{y} \in\Lambda} \bigl(F(x) - F(\overline{x})\bigr)
\sum_{x'\in B} g(x,x') H(x')
\\
& = &\frac{1}{|\Lambda|^2} \sum_{y, \overline{y} \in\Lambda
, x'\in B} \bigl(g(x,x') - g(\overline{x},x')\bigr) F(x) H(x') .
\nonumber
\end{eqnarray}
By~(\ref{3.15}) we thus find that with the notation $\pi_\IZ(x') = z'$,
\begin{eqnarray}\label{3.20}
| \langle F(GH)\rangle|_z &\le&\frac{c}{|\Lambda|^2} \sum
_{y, \overline{y} \in\Lambda, x'\in B} \frac{1}{1 +
|z-z'|^2} \|F\|_\infty \|H\|_\infty\nonumber\\[-8pt]\\[-8pt]
&\le& c \|F\|_\infty \|H\|
_\infty,\nonumber
\end{eqnarray}
and~(\ref{3.18}) follows.
\end{pf}

\section{Occupation times of long rods in ${\IZ^3}$}\label{sec4}

In this section we relate the field of occupation times of long rods in
$\IZ^3$ by random interlacements at a suitably scaled level [see (\ref
{4.1}) below] with the two-dimensional free field pinned at the origin
introduced in~(\ref{1.29}). The main results are stated in 
Theorems~\ref{theo4.2} and~\ref{4.9}. The approach is roughly the following.
By Theorem~\ref{theo2.1} we can express the characteristic functionals
of the scaled fields of occupation times of the long rods as
exponentials of certain power series. The main task is to control the
asymptotic behavior of these power series. This analysis is carried out
in the central Theorem~\ref{theo4.1} as well as in the simpler
Theorem~\ref{theo4.8}. Throughout this section we assume that $d=3$. The
constants depend on the finite subset $\Lambda$ of $\IZ^2$ [cf. (\ref
{3.1}) and above~(\ref{4.3})] as well as on the function $W$ with
support in $\Lambda$ that appears in~(\ref{4.3}). As in Section~\ref{sec3} we
denote by $\|\cdot\|_\infty$ the supremum norm $\| \cdot\|_{L^\infty
(B)}$, with $B$ as in~(\ref{3.1}).

We consider $\alpha> 0$, and a positive sequence $\gamma_N$ tending
to infinity. We will analyze the random fields of occupation times of
the long rods $J_y$, $y \in\IZ^2$, with $J_y = \{y\} \times J = \{y\}
\times\{1,\ldots,N\} \subseteq\IZ^3$, by random interlacements at
the scaled levels
%
%
\begin{equation}\label{4.1}
u_N = \alpha \frac{\log N}{N} \quad\mbox{and}\quad u'_N = \gamma
_N \frac{\log N}{N}\qquad\mbox{with } N > 1.
\end{equation}
The corresponding occupation times of the rods $J_y$, $y \in\IZ^2$, are
%
%
\begin{equation}\label{4.2}
\cL_{y,N} = \sum_{x \in J_y} L_{x,u_N},\qquad \cL'_{y,N} = \sum_{x
\in J_y} L_{x,u'_N}.
\end{equation}
Let us point out that sequences of levels converging faster to zero
than $u_N$ are not interesting in the present context; see Remark \ref
{rem4.3} below.

As in Section~\ref{sec3}, we consider some $\Lambda\subset\subset\IZ^2$
containing $0$. Further we introduce a function $W$ on $\IZ^2$ such that:
%
%
\begin{equation}\label{4.3}
\mbox{(i)\quad} W(y) = 0 \mbox{ outside $\Lambda$};\qquad
\mbox{(ii)\quad} \sum_y W(y) = 0 .
\end{equation}
We define the functions on $\IZ^3$,
\begin{eqnarray}\label{4.4}
V_N(x) &=& \frac{1}{\sqrt{\log N}} W(y) 1_J(z),\nonumber\\[-8pt]\\[-8pt]
V'_N(x) &=& \frac{1}{\sqrt{\gamma_N}} V_N(x)\qquad\mbox{with $x =
(y,z)$},\nonumber
\end{eqnarray}
so that
%
%
\begin{equation}\label{4.5}
\cL_N \stackrel{\mathrm{def}}{=} \sum_{y \in\IZ^2} W(y) \frac
{\cL_{y,N}}{\sqrt{\log N}} = \sum_{x \in\IZ^3} V_N(x) L_{x,u_N}
\end{equation}
and similarly
%
%
\begin{equation} \label{4.6}
\cL'_N \stackrel{\mathrm{def}}{=} \sum_{y \in\IZ^2} W(y) \frac
{\cL'_{y,N}}{\sqrt{\gamma_N \log N}} = \sum_{x \in\IZ^3} V'_N(x)
L_{x,u'_N}.
\end{equation}
It follows from Theorem~\ref{theo2.1} and Remark~\ref{rem2.2} that
%
%
\begin{equation}\label{4.7}
\IE[\exp\{z \cL_N\}] = \exp\biggl\{ \sum_{n \ge1} a_N(n)
z^n\biggr\}\qquad \mbox{for $|z| < r_N$ in $\IC$}
\end{equation}
with $r_N > 0$ and
%
%
\begin{equation}\label{4.8}
a_N(n) = u_N (V_N, (GV_N)^{n-1} 1) \qquad\mbox{for $n \ge1$}.
\end{equation}
As a result of the centering condition~(\ref{4.3})(ii) we have
%
%
\begin{equation}\label{4.9}
a_N(1) = 0.
\end{equation}
In a similar fashion we have
%
%
\begin{equation}\label{4.10}
\IE[\exp\{z \cL'_N\}] = \exp\biggl\{ \sum_{n \ge1} a'_N(n)
z^n\biggr\}\qquad \mbox{for $|z| < r'_N$ in $\IC$}
\end{equation}
with $r'_N > 0$ and
%
%
\begin{equation}\label{4.11}
a'_N(n) = u'_N (V'_N, (GV'_N)^{n-1} 1)
\stackrel{\mbox{\fontsize{8.36}{8.36}\selectfont{(\ref{4.1}),
(\ref{4.4})}}}{=} \frac{1}{\alpha} \gamma_N^{1 - {n/2}}
a_N(n)\qquad\mbox{for $n \ge1$}.\hspace*{-35pt}
\end{equation}
The heart of the matter for the proof of Theorem~\ref{theo4.2} lies in
the analysis of the large $N$ behavior of the coefficients $a_N(n)$, $n
\ge1$. There is a dichotomy between the case of odd $n$, with an
asymptotic vanishing of $a_N(n)$, and even~$n$, with a positive limit
of $a_N(n)$, as $N$ goes to infinity. The crucial controls are
contained in the next theorem. We recall the convention concerning
constants stated at the beginning of this section.
\begin{theorem}\label{theo4.1}
%
%
\begin{eqnarray}
\label{4.12}
&\displaystyle  |a_N(n)| \le\alpha c_0^n\qquad \mbox{for all } n \ge1, N> 1,&
\\
\label{4.13}
&\mbox{for any }k \ge0\qquad\displaystyle \lim_N a_N(2k+1) =0,&
\\
\label{4.14}
&\mbox{for any } k \ge1\qquad\displaystyle  \lim_N a_N(2k) =\frac
{\alpha}{2} \cE(W) \biggl(\frac{3}{2 \pi} \cE(W)
\biggr)^{k-1},&
\end{eqnarray}
where we have set [see~(\ref{1.6}) for notation]
%
%
\begin{equation}\label{4.15}
\cE(W) = - 3 \sum_{y,y'} W(y) W(y') a(y' - y).
\end{equation}
\end{theorem}

Note that due to~(\ref{4.3})(ii) we can express $\cE(W)$ in
terms of the two-dimensional Gaussian free field $\psi_y$, $y \in\IZ
^2$, introduced in~(\ref{1.29}), via the formula
%
%
\begin{equation}\label{4.16}
\cE(W) = E\biggl[\biggl(\sum_{y \in\IZ^2} W(y) \psi_y\biggr)^2\biggr].
\end{equation}
Before turning to the proof of Theorem~\ref{theo4.1}, we first explain
how this theorem enables us to derive the convergence in law of the
appropriately\vadjust{\goodbreak} scaled fields $\cL_{y,N} - \cL_{0,N}$, $y \in\IZ^2$
and $\cL'_{y,N} - \cL'_{0,N}$, $y \in\IZ^2$. We tacitly endow $\IR
^{\IZ^2}$ with the product topology, so that the convergence stated in
Theorem~\ref{theo4.2} actually corresponds to the convergence in
distribution of all finite-dimensional marginals of the relevant random
fields. The main result of this section is the next theorem, which
proves~(\ref{0.8}) and~(\ref{0.10}).
\begin{theorem}\label{theo4.2}
%
%
\begin{equation}\label{4.17}
\begin{tabular}{p{315pt}}
As $N$ goes to infinity, $(\frac{\cL_{y,N} - \cL
_{0,N}}{\sqrt{\log N}})_{y \in\IZ^2}$, converges in
distribution to the random field $(R \psi_y)_{y \in\IZ^2}$,
\end{tabular}\hspace*{-37pt}
\end{equation}
where $R$ and $(\psi_y)_{y \in\IZ^2}$ are independent and
%
%
\begin{eqnarray}
\label{4.18}
&& \mbox{$R$ is $\mathrm{BES}^0(\sqrt{\alpha}, \frac{3}{2 \pi}
)$-distributed,}
\\
\label{4.19}
&&\mbox{$(\psi_y)_{y \in\IZ^2}$ is the centered Gaussian field
introduced in~(\ref{1.29}).}
\end{eqnarray}
Moreover,
%
%
\begin{equation}\label{4.20}
\begin{tabular}{p{300pt}}
as $N$ goes to infinity, $(\frac{\cL'_{y,N} - \cL
'_{0,N}}{\sqrt{Nu'_N}})_{y \in\IZ^2}$ converges in
distribution to $(\psi_y)_{y \in\IZ^2}$.
\end{tabular}\hspace*{-35pt}
\end{equation}
\end{theorem}
\begin{pf*}{Proof 
(\textup{assuming Theorem~\ref{theo4.1})}}
We begin with the proof of~(\ref{4.17}). We consider $W$
as in~(\ref{4.3}) and $\cL_N$ as in~(\ref{4.5}). By~(\ref{4.7}),
(\ref{4.12}) and Remark~\ref{rem2.2}, we see that $\exp\{z \cL_N\}$
is integrable for any $z$ with $c_0 |\operatorname{Re} z| < 1$, and that
%
%
\begin{equation}\label{4.21}\quad
\IE[\exp\{ z \cL_N\}] = \exp\biggl\{\sum_{n \ge1} a_N(n) z^n\biggr\}
\qquad \mbox{for any $|z| < c_0^{-1}$ in $\IC$}.
\end{equation}
In particular~(\ref{4.12}) implies that
%
%
\begin{equation}\label{4.22}
\sup_N \IE[\operatorname{cosh} (r \cL_N)] < \infty\qquad \mbox{when
$r < c_0^{-1}$}.
\end{equation}
As a result the laws of the variables $\cL_N$ are tight, and the
variables $\exp\{z \cL_N\}$, \mbox{$N > 1$}, with $|{\operatorname{Re}
z}| \le r < c_0^{-1}$, are uniformly integrable. If along some
subsequence~$N_k$, $k \ge1$, the variables $\cL_{N_k}$ converge in
distribution to $\cL $, it follows from Theorem~5.4, page 32 in
\cite{Bill68}, that for $|z| < c_0^{-1}$,
%
%
\begin{eqnarray}\label{4.23}
E[\exp\{z \cL\}] & = &\lim_k \IE[\exp\{z \cL_{N_k}\}]
\nonumber\\
& = &\lim_k \exp\biggl\{\sum_{n \ge1} a_{N_k} (n) z^n
\biggr\}
\\
& = &\exp\biggl\{\frac{\alpha}{2} \cE(W) z^2 \Big/ \biggl(1 -
\frac{3}{2 \pi} \cE(W)z^2\biggr)\biggr\}
\nonumber
\end{eqnarray}
using Theorem~\ref{theo4.1} in the last equality. This determines the
characteristic function of the law of $\cL$, and by~(\ref{1.30})
shows that $\cS$ has same distribution as $R\psi$ where $R,\psi$ are
independent variables with $R \,\mathrm{BES}^0(\sqrt{\alpha}, \frac{3}{2\pi
})$-distributed,\vadjust{\goodbreak} and $\psi$ a centered Gaussian variable with zero
mean and variance $\cE(W)$. This proves that for any $W$ as in (\ref
{4.3}), $\cL_N$ converges in distribution to $R \psi$ as above, when
$N$ tends to infinity. In view of~(\ref{4.16}), this completes the
proof of~(\ref{4.17}).

The proof of~(\ref{4.20}) is analogous. Due to~(\ref{4.11}), we know
that $a'_N(n) = \frac{1}{\alpha} \gamma_N^{1 -
{n/2}}a_N(n)$, and in particular\vspace*{1pt} $a'_N(2) = \frac{1}{\alpha}
a_N(2)$ converges to $\frac{1}{2} \cE(W)$,\break whereas for $n\not=2$,
$a'_N(n)$ converges to zero as $N$ goes to infinity. We can\vspace*{1pt} use the
same arguments as above and find that $\cL'_N$ converges in
distribution to a variable $\cL'$ such that for small $|z|$ in $\IC$,
$\IE[\exp\{z \cL'\}] = \exp\{\frac{1}{2} \cE(W) z^2\}$. The
claim~(\ref{4.20}) then follows immediately.
\end{pf*}
\begin{remark}\label{rem4.3}
For the kind of limit theorems discussed here, sequences of levels
converging to zero faster than $u_N$ lead to trivial results, as we now
explain. In a standard way (see, e.g., Remark 3.1(3) in
\cite{SidoSzni09b}), one has the bound $\operatorname{cap}(J_y) \le c \frac
{N}{\log N}$, for all $y \in\IZ^2$. If we pick $u''_N$ so that $u''_N
= o(\frac{\log N}{N})$, then~(\ref{2.20}) implies that for all $y \in
\IZ^2$, with probability tending to $1$ as $N$ goes to infinity, the
interlacement at level $u''_N$ does not intersect $J_y$. In particular,
if we define $\cL''_{y,N}$ in analogy to $\cL_{y,N}$ in~(\ref{4.2})
with $u''_N$ in place of $u_N$, the random field $(\cL''_{y,N})_{y \in
\IZ^2}$ converges\vspace*{1pt} in distribution to the constant field equal to zero,
as $N$ goes to infinity.
\end{remark}
\begin{pf*}{Proof of Theorem~\ref{theo4.1}} We recall the convention
concerning constants and the notation \mbox{$\| \cdot\|_\infty$} stated at
the beginning of this section. The linear operators under consideration
throughout the proof will be restricted to the space of functions
vanishing outside $B$. The centering condition~(\ref{4.3})(ii) and
the ensuing identity $\langle V_N\rangle= 0$ play a crucial role. We
will first prove~(\ref{4.12}), (\ref{4.13}). Our first step is to
control the norm of the operators $(GV_N)^2$.
\begin{lemma}\label{lem4.4}
%
%
\begin{equation}\label{4.24}
\|(GV_N)^2 \|_{L^\infty(B) \rightarrow L^\infty(B)} \le c_1 .
\end{equation}
\end{lemma}
\begin{pf}
Given a function $F$ on $B$, we write in the notation of~(\ref{3.2}),~(\ref{3.3})
\begin{eqnarray}\label{4.25}
(GV_N)^2 F &=& A_1 + A_2 + A_3\qquad \mbox{where}
\nonumber\\
A_1 &=& GV_N G(V_N F - \langle V_N F \rangle),
\nonumber\\[-8pt]\\[-8pt]
A_2 &=& GV_N [G\langle V_N F \rangle]_0,
\nonumber\\
A_3 &=& GV_N (G \langle V_N F\rangle- [G \langle V_N
F\rangle]_0).\nonumber
\end{eqnarray}
With the help of Lemma~\ref{lem3.1} we see that
%
%
\begin{equation}\label{4.26}\qquad
\|A_1\|_\infty\stackrel{\mbox{\fontsize{8.36}{8.36}\selectfont{(\ref{3.8}),~(\ref{4.4})}}}{\le} c \sqrt
{\log N} \|G(V_N F - \langle V_N F\rangle)\|_\infty
\stackrel{\mbox{\fontsize{8.36}{8.36}\selectfont{
{(\ref{3.10}),~(\ref{4.4})}}}}{\le} c \|F\|_\infty.
\end{equation}
Since $\langle V_N[G\langle V_N F\rangle]_0\rangle= 0$, we also
find that
%
%
\begin{equation}\label{4.27}
\|A_2\|_\infty\stackrel{\mbox{\fontsize{8.36}{8.36}\selectfont{(\ref{3.10}),
(\ref{4.4})}}}{\le}
\frac{c}{\sqrt{\log N}} \|G \langle V_N F\rangle\|_\infty
\stackrel{\mbox{\fontsize{8.36}{8.36}\selectfont{(\ref{3.8}),~(\ref{4.4})}}}{\le} c \|F\|_\infty.
\end{equation}
Finally we have
\begin{eqnarray}\label{4.28}
\|A_3\|_\infty
&\stackrel{\mbox{\fontsize{8.36}{8.36}\selectfont{(\ref{3.8}),~(\ref{4.4})}}}
{\le}& c \sqrt
{\log N} \|G(V_N F - [G \langle V_N F\rangle]_0)
\|_\infty\nonumber\\[-8pt]\\[-8pt]
&\stackrel{\mbox{\fontsize{8.36}{8.36}\selectfont{(\ref{3.12}),
(\ref{4.4})}}}{\le}& c \|F\|_\infty.\nonumber
\end{eqnarray}
Collecting~(\ref{4.25})--(\ref{4.28}), the claim~(\ref{4.24}) now follows.
\end{pf}

Before proving~(\ref{4.12}),~(\ref{4.13}) we still need the following
lemma, which shows that the kernel of the linear operator $F
\rightarrow\langle V_N F\rangle$, is almost invariant under
$(GV_N)^2$. We will later see [cf.~(\ref{4.61})] that the function $F=
1_B$, which\vspace*{1pt} belongs to this kernel, is in an appropriate sense, close
to being an eigenvector of $(GV_N)^2$.
\begin{lemma}
When $F$ is a function on $B$, one has
%
%
\begin{equation}\label{4.29}
\|\langle V_N(GV_N)^2 F\rangle\|_\infty\le c_2 \| \langle V_N
F\rangle\|_\infty+ \frac{c_3}{(\log N)^{{3/2}}} \|F\|
_\infty.
\end{equation}
\end{lemma}
\begin{pf}
We use~(\ref{4.25}) and write
%
%
\begin{equation}\label{4.30}
\langle V_N(GV_N)^2 F\rangle= \langle V_N A_1 \rangle+ \langle V_N
A_2 \rangle+ \langle V_N A_3 \rangle.
\end{equation}
With the help of Lemma~\ref{lem3.2} we find that
%
%
\begin{eqnarray}
\label{4.31}
\| \langle V_N A_1 \rangle\|_\infty& \stackrel{\mbox{\fontsize{8.36}{8.36}
\selectfont{(\ref{3.18})}}}{\le}&
\frac{c}{\sqrt{\log N}} \|V_N G(V_N F - \langle V_N F\rangle
)\|_\infty\nonumber\\[-8pt]\\[-8pt]
&\stackrel{\mbox{\fontsize{8.36}{8.36}\selectfont{(\ref{3.10})}}}{\le}&
\frac{c}{(\log N)^{{3/2}}} \|F\|_\infty,\nonumber
\\
\label{4.32}
\| \langle V_N A_2 \rangle\|_\infty& \stackrel{\mbox{\fontsize{8.36}{8.36}
\selectfont{(\ref{3.18})}}}{\le}&
\frac{c}{\sqrt{\log N}} \|V_N [G\langle V_N F \rangle]_0 \|
_\infty\stackrel{\mbox{\fontsize{8.36}{8.36}\selectfont{(\ref{3.8})}}}{\le} c \|\langle V_N F\rangle\|
_\infty
\end{eqnarray}
and that
\begin{eqnarray}\label{4.33}
\| \langle V_N A_3 \rangle\|_\infty& \stackrel{\mbox{\fontsize{8.36}{8.36}
\selectfont{(\ref{3.18})}}}{\le}&
\frac{c}{\sqrt{\log N}} \|V_N (G\langle V_N F\rangle- [G
\langle V_N F \rangle]_0)\|_\infty
\nonumber\\[-8pt]\\[-8pt]
& \stackrel{\mbox{\fontsize{8.36}{8.36}\selectfont{(\ref{3.12})}}}{\le}&
\frac{c}{(\log N)^{{3/2}}} \| F \|_\infty.\nonumber
\end{eqnarray}

Collecting~(\ref{4.30})--(\ref{4.33}), we obtain~(\ref{4.29}).
\end{pf}

We now prove~(\ref{4.12}),~(\ref{4.13}). As a result of (\ref
{4.24}),~(\ref{4.29}) we see that for $k \ge1$ and $F$ a function on
$B$, one has
\begin{eqnarray}\label{4.34}
\|\langle V_N (G V_N)^{2k} F \rangle\|_\infty
&\stackrel{\mbox{\fontsize{8.36}{8.36}
\selectfont{(\ref{4.29}),~(\ref{4.24})}}}{\le}&
c_2 \bigl\| \bigl\langle V_N(GV_N)^{2(k-1)} F
\bigr\rangle\bigr\|_\infty\nonumber\\[-8pt]\\[-8pt]
&&{}+ \frac{c_3}{(\log N)^{{3/2}}}
c_1^{k-1} \|F\|_\infty\nonumber
\end{eqnarray}
and by induction
\begin{eqnarray*}
& \le &c_2^k \| \langle V_N F\rangle\|_\infty+ \frac{c_3}{(\log
N)^{3/2}} (c^{k-1}_1 + c_2 c_1^{k-2} +\cdots+ c_2^{k-1})
\|F\|_\infty
\\
& \le &c_2^k \| \langle V_N F\rangle\|_\infty+ \frac{c^k}{(\log
N)^{3/2}} \|F\|_\infty.
\end{eqnarray*}
Keeping in mind that $\langle V_N\rangle= 0$, we thus see that for $k
\ge1$,
%
%
\begin{eqnarray}\label{4.35}
|a_N(2k+1)|& \stackrel{\mbox{\fontsize{8.36}{8.36}\selectfont{(\ref{4.1}),
(\ref{4.8})}}}{=} & \frac
{\alpha\log N}{N} |(V_N, (GV_N)^{2k} 1)| \nonumber\\
&\le&\alpha\log N |\Lambda
| \|\langle V_N (GV_N)^{2k} 1 \rangle\|_\infty
\\
& \stackrel{\mbox{\fontsize{8.36}{8.36}\selectfont{(\ref{4.34})}}}{\le}
& \frac{\alpha c^k |\Lambda |}{\sqrt{\log N}} . \nonumber
\end{eqnarray}
Together with~(\ref{4.9}), this proves~(\ref{4.13}) as well as (\ref
{4.12}) for odd $n$. When $n = 2k$, with $k \ge1$, we note that
%
%
\begin{eqnarray}\label{4.36}\qquad
|a_N(2k)| & = & \frac{\alpha\log N}{N} \bigl|\bigl(V_N,
(GV_N)(GV_N)^{2(k-1)} 1\bigr)\bigr|
\nonumber\\
& = & \frac{\alpha\log N}{N} \bigl| \bigl(GV_N, V_N (GV_N)^{2(k-1)} 1\bigr)\bigr|
\qquad\mbox{(by symmetry of $G$)}
\\
& \le & \alpha\log N |\Lambda| \|GV_N \|_\infty\bigl\|
V_N(GV_N)^{2(k-1)} 1\bigr\|_\infty\stackrel{\mbox{\fontsize{8.36}{8.36}
\selectfont{(\ref{3.10}),~(\ref{4.24})}}}{\le} \alpha c^k .
\nonumber
\end{eqnarray}
The proof of~(\ref{4.12}) is now complete.

There remains to prove~(\ref{4.14}), that is, to analyze the large $N$
behavior of the even coefficients $a_N(2k)$. To motivate the next lemma
we recall that
%
%
\begin{equation}\label{4.37}
a_N(2) \stackrel{\mbox{\fontsize{8.36}{8.36}\selectfont{(\ref{4.8})}}}{=} \frac{\alpha\log N}{N}
(V_N, GV_N) \stackrel{\mbox{\fontsize{8.36}{8.36}\selectfont{(\ref{3.2})}}}{=} \frac{\alpha\log N}{N}
| \Lambda| \sum_{z \in J} \langle V_N GV_N\rangle_z .
\end{equation}

\begin{lemma}\label{lem4.6}
There exists a function $\Gamma(\cdot)$ on $\IN$ tending to $0$ at
infinity such that
%
%
\begin{equation}\label{4.38}
\log N \langle V_N G V_N\rangle_z = \fr \frac{\cE
(W)}{|\Lambda|} + f_N(z)\qquad \mbox{for $z \in J$}
\end{equation}
with $|f_N(z)| \le\Gamma(d(z,J^c))$.

Moreover if one defines
%
%
\begin{equation}\label{4.39}
\tau_N = \frac{1}{2 \log N} \sum_{|z| \le N} g((0,z))\qquad
\mbox{so that }\lim_N \tau_N = \frac{3}{2 \pi}\mbox{ by
(\ref{1.9})},
\end{equation}
one has the identity on $B$
%
%
\begin{equation}\label{4.40}
(GV_N)^2 1_B - \tau_N \cE(W) 1_B = \frac{1}{\log N} G(f_N
\circ\pi_\IZ) + k_N,
\end{equation}
where for $x = (y,z) \in B$,
%
%
\begin{equation}\label{4.41}
|k_N|(x) \le\frac{c}{\log N} \log\bigl(N / d(z,J^c)\bigr).
\end{equation}
\end{lemma}
\begin{pf}
We begin with the proof of~(\ref{4.38}). We note that for $z \in J$,
$x = (y,z)$ in $B$ and $x'= (y',z')$ in $B$, one has
%
%
\begin{eqnarray}\label{4.42}
\log N \langle V_N GV_N\rangle_z &=& \frac{1}{|\Lambda|} \sum
_{y \in\Lambda, x' \in B} W(y) g(x,x') W(y')
\nonumber\\
&\stackrel{\mbox{\fontsize{8.36}{8.36}\selectfont{(\ref{4.3})(ii)}}}{=}
&\frac{1}{|\Lambda|} \sum_{y,y' \in\Lambda, z' \in J} W(y) \bigl(g(x,x') -
g(x,(y,z'))\bigr) W(y')
\nonumber\\
&=& \frac{1}{|\Lambda|} \sum_{y,y' \in\Lambda} W(y) W(y')
\sum_{z' \in J} \bigl(g(x'-x) - g\bigl((0,z'-z)\bigr)\bigr)
\\
&\stackrel{\mbox{\fontsize{8.36}{8.36}\selectfont{(\ref{1.10})}}}{=} & \frac{1}{|\Lambda|} \sum
_{y,y' \in\Lambda} W(y) W(y') \biggl(b_N(y' - y,z) - \frac
{3}{2} a(y' - y)\biggr)
\nonumber\\
&\stackrel{\mbox{\fontsize{8.36}{8.36}\selectfont{(\ref{4.15})}}}{=} & \fr\frac{\cE(W)}{| \Lambda|}
+ f_N(z),\nonumber
\end{eqnarray}
where we have set
\[
f_N(z) = \frac{1}{|\Lambda|} \sum_{y,y' \in\Lambda} W(y)
W(y') b_N(y' - y, z)\qquad \mbox{for $z \in J$}.
\]
The estimate in the second line of~(\ref{4.38}) is now an immediate
consequence of~(\ref{1.11}). This completes the proof of~(\ref{4.38}).

We then turn to the proof of~(\ref{4.40}),~(\ref{4.41}). We write
%
%
\begin{equation}\label{4.43}\qquad
(GV_N)^2 1_B = G \langle V_N GV_N\rangle+ G(V_N GV_N - \langle
V_N GV_N\rangle) \stackrel{\mathrm{def}}{=} a_1 + a_2 .
\end{equation}
We know that
%
%
\begin{equation}\label{4.44}
\|a_2\|_\infty\stackrel{\mbox{\fontsize{8.36}{8.36}\selectfont{(\ref{3.10})}}}{\le} c \|V_N GV_N\|
_\infty\stackrel{\mbox{\fontsize{8.36}{8.36}\selectfont{(\ref{3.10})}}}{\le} c' / \log N,
\end{equation}
and by~(\ref{4.38}) we find that
\begin{eqnarray}\label{4.45}
a_1 & = & \fr \frac{\cE(W)}{| \Lambda|} \frac
{G1_B}{\log N} + \frac{1}{\log N} G(f_N \circ\pi_\IZ)
\nonumber\\[-8pt]\\[-8pt]
& = & \tau_N \cE(W) 1_B + \frac{1}{\log N} G(f_N \circ\pi
_\IZ) + r_N,
\nonumber
\end{eqnarray}
where for $x = (y,z) \in B$ we have set
%
%
\begin{equation}\label{4.46}
r_N(x) = \fr \frac{\cE(W)}{| \Lambda|} \frac{G1
_B}{\log N} - \fr \frac{\cE(W)}{\log N} \sum_{|z' - z|
\le N} g\bigl((0,z' - z)\bigr).
\end{equation}
We thus see that (recall $J_y = \{y\} \times J)$
%
%
\begin{eqnarray}\label{4.47}
&&|r_N(x)| \le \frac{c}{\log N} \biggl\| G \biggl(\frac{1_B}{|
\Lambda|} - 1_{J_0}\biggr)\biggr\|_\infty+ \frac{c}{\log N} \|G
1_{J_0} - [G1_{J_0}]_0\|_\infty
\nonumber\\
&&\hphantom{|r_N(x)| \le}{}
+ \frac{c}{\log N} \sum_{|z' - z| \le N, z' \notin J}
g\bigl((0,z' - z)\bigr) \\
&&\hspace*{6.8pt}\stackrel{\mbox{\fontsize{8.36}{8.36}\selectfont{(\ref{3.10}),
(\ref{3.12}),~(\ref{1.3})}}}{\le}
\frac{c}{\log N} \log\bigl(N/ d(z,J^c)\bigr) .
\nonumber
\end{eqnarray}
Collecting~(\ref{4.43})--(\ref{4.47}) we have completed the proof of
(\ref{4.40}),~(\ref{4.41}).
\end{pf}

As a result of~(\ref{4.37}),~(\ref{4.38}) we see that
%
%
\begin{equation}\label{4.48}
a_N(2) = \frac{\alpha}{2} \cE(W) +
\frac{\alpha|\Lambda|}{N} \sum_{z \in J} f_N(z)
\ulN\frac{\alpha}{2} \cE(W).
\end{equation}
This proves~(\ref{4.14}) in the case $k=1$. To handle the case $k >
1$, we will need to control the propagation of boundary effects
corresponding to terms with $\IZ$-component close to the complement of
$J$, when proving the convergence of $a_N(2k) = \alpha\frac{\log
N}{N} (V_N,(GV_N)^{2k-1} 1)$ for $N \rightarrow\infty$. The next lemma will
be useful for this purpose. We first introduce some notation.

We consider nonempty sub-intervals of $J$,
%
%
\begin{equation}\label{4.49}
I_0 \subseteq I_1 \subseteq I_2 \subsetneq J
\end{equation}
and define
%
%
\begin{equation}\label{4.50}
L_0 = d(I_0, J \setminus I_1) \ge1,\qquad L_1 = d(I_1, J \setminus I_2)
\ge1,
\end{equation}
as well as
%
%
\begin{equation}\label{4.51}
C_0 = \Lambda\times I_0 \subseteq C_1 = \Lambda\times I_1 = C_2 =
\Lambda\times I_2 \subseteq B.
\end{equation}
We have the following variation on Lemma~\ref{lem4.4}.
\begin{lemma}\label{lem4.7}
For $F$ a function on $B$, one has
\begin{eqnarray}\label{4.52}\qquad
&&
\|1_{C_0} (GV_N)^2 F\|_\infty\nonumber\\[-9pt]\\[-9pt]
&&\qquad\le c_4 \|1_{C_2} F\|_\infty+ c_5
\biggl(\frac{1}{L_0} + \frac{1}{L_1} + \frac{1}{\log N} \log
\biggl(\frac{(N+1)^2}{L_0 L_1}\biggr)\biggr) \| F\|_\infty.\nonumber
\end{eqnarray}
\end{lemma}
\begin{pf}
By Lemma~\ref{lem4.4} we can assume that $F$ vanishes on $C_2$. We use
the decomposition~(\ref{4.25}) of $(GV_N)^2 F$. We find that
\begin{eqnarray}\label{4.53}
&&\|1_{C_0} A_1\|_\infty = \|1_{C_0} GV_N G(V_N F - \langle V_N F\rangle
)\|_\infty
\nonumber\\[-2pt]
&&\hphantom{\|1_{C_0} A_1\|_\infty }\le\|1_{C_0} GV_N 1_{C_1}
G(V_N F - \langle V_N F\rangle)\|_\infty\nonumber\\[-2pt]
&&\hphantom{\|1_{C_0} A_1\|_\infty \le}
{} + \| 1_{C_0} GV_N 1_{B \setminus C_1} G(V_N F - \langle V_N F\rangle
)\|_\infty
\nonumber\\[-9pt]\\[-9pt]
&&\hspace*{36pt}\stackrel{\mbox{\fontsize{8.36}{8.36}\selectfont{(\ref{3.8}),~(\ref{3.9})}}}
{\le} c \sqrt{\log N} \|
1_{C_1} G(V_N F - \langle V_N F\rangle)\|_\infty
\nonumber\\[-2pt]
&&\hphantom{\|1_{C_0} A_1\|_\infty \le}+ \frac{c}{\sqrt{\log N}} \log\biggl(\frac
{N+1}{L_0}\biggr) \|G(V_N F- \langle V_N F\rangle)\|_\infty
\nonumber\\[-2pt]
&&\hspace*{32pt}\stackrel{\mbox{\fontsize{8.36}{8.36}\selectfont{(\ref{3.11}),
(\ref{3.10})}}}{\le} \frac{c}{L_1} \|
F\|_\infty+ \frac{c}{\log N} \biggl(\frac{N+1}{L_0}
\biggr) \| F\|_\infty
\nonumber
\end{eqnarray}
using the fact that $F$ vanishes on $C_2$ in the last step. In a
similar fashion we find that
\begin{eqnarray}\label{4.54}
&&\|1_{C_0} A_2\|_\infty= \|1_{C_0} GV_N [G\langle V_N F \rangle]_0 \|
_\infty\nonumber\\[-2pt]
&&\hspace*{32.36pt}\stackrel{\mbox{\fontsize{8.36}{8.36}\selectfont{(\ref{3.10}),
(\ref{3.11})}}}{\le}
\frac{c}{\sqrt{\log N}} \|1_{C_1} [G \langle V_N F\rangle
]_0 \|_\infty\nonumber\\[-9pt]\\[-9pt]
&&\hphantom{\|1_{C_0} A_2\|_\infty=}{}+ \frac{c}{\sqrt{\log N}} \frac{1}{L_0}
\|1_{B \setminus C_1} [G\langle V_N F\rangle]_0 \|_\infty
\nonumber\\[-2pt]
&&\hspace*{36pt}\stackrel{\mbox{\fontsize{8.36}{8.36}\selectfont{(\ref{3.9}),
(\ref{3.8})}}}{\le}
\frac{c}{\log N} \log\biggl(\frac{N+1}{L_1}\biggr) \| F\|
_\infty+ \frac{c}{L_0} \| F\|_\infty,
\nonumber
\end{eqnarray}
where once again we have used that $F$ vanishes on $C_2$ in the last
step. Finally we have
\begin{eqnarray}\label{4.55}
&&\|1_{C_0} A_3\|_\infty= \|1_{C_0} GV_N (G\langle V_N F \rangle-
[G\langle V_N F\rangle]_0)\|_\infty\nonumber\\[-2pt]
&&\hspace*{36.49pt}\stackrel{\mbox{\fontsize{8.36}{8.36}\selectfont{(\ref{3.8}), (\ref
{3.9})}}}{\le}
c \sqrt{\log N} \|1_{C_1} (G \langle V_N F\rangle- [G \langle V_N
F\rangle]_0)\|_\infty
\nonumber\\[-9pt]\\[-9pt]
&&\hphantom{\|1_{C_0} A_3\|_\infty=}{}+ \frac{c}{\sqrt{\log N}}
\log\biggl(\frac{N+1}{L_0}\biggr) \|1_{B \setminus C_1} (G \langle V_N
F\rangle- [G \langle V_NF\rangle]_0)\|_\infty
\nonumber\\[-2pt]
&&\hspace*{32pt}\stackrel{\mbox{\fontsize{8.36}{8.36}\selectfont{(\ref{3.13}),
(\ref{3.12})}}}{\le} \frac{c}{L_1} \|F\| _\infty+ \frac{c}{\log N}
\log\biggl(\frac{N+1}{L_0}\biggr) \|F\|_\infty \nonumber
\end{eqnarray}
using that $F$ vanishes on $C_2$ in the last step.

Collecting~(\ref{4.53})--(\ref{4.55}), we obtain~(\ref{4.52}).\vadjust{\goodbreak}
\end{pf}

We now introduce the following sequence of possibly empty sub-intervals
of $J$:
%
%
\begin{eqnarray}\label{4.56}
J_k = \bigl\{i \ge1; 1 + 4k\bigl[N e^{-\sqrt{\log N}}\bigr] \le i \le N - 4k \bigl[N
e^{-\sqrt{\log N}}\bigr]\bigr\}\nonumber\\[-8pt]\\[-8pt]
&&\eqntext{\mbox{for $k \ge0$}.}
\end{eqnarray}
There is some freedom in the above definition. The proof below would
work with minor changes if one replaces $N e^{-\sqrt{\log N}}$ by
$N^{1-\varepsilon_N}$, with $\varepsilon_N \rightarrow0$, and
$\varepsilon_N
\log N \rightarrow\infty$.

Setting $B_k = \Lambda\times J_k$, we find that for any $k \ge1$,
when $N \ge c(k)$,
\[
\varnothing\not= B_k \subsetneq B_{k-1} \subsetneq\cdots\subsetneq B_1
\subsetneq B
\]
and that with the notation~(\ref{4.39})
\begin{eqnarray}\label{4.57}\quad
&&\bigl\|1_{B_k} \bigl((GV_N)^{2k} 1_B - (\tau_N \cE(W))^k 1_B\bigr)\bigr\|_\infty\nonumber\\
&&\qquad\le
\sum^{k-1}_{m=0} \bigl\|1_{B_k} \bigl((GV_N)^{2(m+1)} (\tau_N \cE
(W))^{k-(m+1)} 1_B \nonumber\\[-8pt]\\[-8pt]
&&\hspace*{90.6pt}{}- (GV_N)^{2m} (\tau_N \cE(W))^{k-m} 1_B\bigr)\bigr\|_\infty
\nonumber\\
&&\qquad\le
\sum^{k-1}_{m=0} (\tau_N \cE(W))^{k-(m+1)} \bigl\|1_{B_{m+1}}(GV)^{2m}
\bigl((G V_N)^2 1_B - \tau_N \cE(W) 1_B\bigr)\bigr\|_\infty.
\nonumber
\end{eqnarray}
We set\vspace*{1pt} $F_N = (GV_N)^2 1_B - \tau_N \cE(W) 1_B$, and now want to
bound $\|1_{B_{m+1}}(GV_N)^{2m}\times F_N\|_\infty$ with the help of (\ref
{4.52}), when $0 < m < k$. To this end we introduce $\wh{J}_m$ with a
similar definition as in~(\ref{4.56}), simply replacing $4k$ by
$2m+2$, so that $J_{m+1} \subseteq\wh{J}_m \subseteq J_m$ play the
role of $I_0 \subseteq I_1 \subseteq I_2$ in~(\ref{4.49}). We note
that $d(J_{m+1}, J \setminus\wh{J}_m)$ and $d(\wh{J}_m, J
\setminus J_m)$ are bigger than $c N e^{-\sqrt{\log N}}$. As a
result the expression inside the parenthesis after $c_5$ in (\ref
{4.52}) is smaller than
\[
c N^{-1} e^{\sqrt{\log N}} + \frac{1}{\log N} \log
\biggl(\frac{(N+1)^2}{N^2} e^{2 \sqrt{\log N}}\biggr) \le\frac
{c}{\sqrt{\log N}}.
\]
It thus follows that for $0 < m < k$
\begin{eqnarray}\label{4.58}
\|1_{B_{m+1}}(GV_N)^{2m} F_N \|_\infty&\stackrel{\mbox{\fontsize{8.36}{8.36}
\selectfont{(\ref{4.52}), (\ref
{4.24})}}}{\le}& c_4 \bigl\|1_{B_m} (GV_N)^{2(m-1)} F_N \bigr\|_\infty\nonumber\\[-8pt]\\[-8pt]
&&{}+ \frac
{c_1^{m-1}}{\sqrt{\log N}} c_6 \|F_N\|_\infty\nonumber
\end{eqnarray}
and by induction
\begin{eqnarray*}
& \le & c_4^m \|1_{B_1} F_N\|_\infty+ \frac{c_6}{\sqrt{\log N}}
(c_1^{m-1} + c_4 c_1^{m-2} +\cdots+ c_4^{m-1}) \|F_N\|_\infty
\\
& \le & c_4^m \|1_{B_1} F_N\|_\infty+ \frac{c^m}{\sqrt{\log N}}
\|F_N\|_\infty.
\end{eqnarray*}
Coming back to the last line of~(\ref{4.57}), we see that for $k \ge
1$, $N \ge c(k)$, each term under the sum, thanks to the above bound
and~(\ref{4.39}), is smaller than $c(k) (\|1_{B_1} F_N\|_\infty+
\frac{1}{\sqrt{\log N}} \| F_N\|_\infty)$. Hence we see that for $k
\ge1$ and $N \ge c(k)$,
%
%
\begin{equation}\label{4.59}
\bigl\|1_{B_k}\bigl((GV_N)^{2k} 1_B - (\tau_N \cE(W))^k 1_B\bigr)\bigr\|_\infty
\le c(k) \biggl(\!\|1_{B_1} F_N\|_\infty+ \frac{1}{\sqrt{\log
N}}\!\biggr),\hspace*{-35pt}
\end{equation}
where we have used the bound $\|F_N\|_\infty\le c$, which follows from
(\ref{4.24}) and (\ref{4.39}). Note that the form of the correction
term $(\log N)^{-1/2}$ in~(\ref{4.59}) mainly reflects our
choice for the intervals $J_k$ in~(\ref{4.56}). In view of Lemma \ref
{lem4.6}, we also find that
\begin{eqnarray}\label{4.60}
\|1_{B_1} F_N \|_\infty
&\stackrel{\mbox{\fontsize{8.36}{8.36}\selectfont{(\ref{4.40})}}}{\le}
&\biggl\|1_{B_1}
\frac{G(f_N \circ\pi_\IZ)}{\log N} \biggr\|_\infty\nonumber\\[-8pt]\\[-8pt]
&&{} + \|1_{B_1}
k_N\|_\infty\longrightarrow0\qquad\mbox{as $N
\rightarrow\infty$},\nonumber
\end{eqnarray}
where we have used the bounds in the second line of~(\ref{4.38}) and
(\ref{4.41}) to conclude in the last step. Since $\tau_N$ converges
to $\frac{3}{2 \pi}$ [see~(\ref{4.39})], we can infer from (\ref
{4.59}),~(\ref{4.60}) that
%
%
\begin{eqnarray}\label{4.61}
\Delta_{k,N} \stackrel{\mathrm{def}}{=} \biggl\|1_{B_k} \biggl((GV_N)^{2k}
1_B - \biggl(\frac{3}{2 \pi} \cE(W)\biggr)^k 1_B\biggr)\biggr\|_\infty
\ulN0\nonumber\\[-8pt]\\[-8pt]
&&\eqntext{\mbox{for each $k \ge1$}.}
\end{eqnarray}
We can now use this estimate to study the asymptotic behavior of
$a_N(2(k+1))$ as $N$ goes to infinity. Indeed one has
%
%
\begin{eqnarray}\label{4.62}
&&\biggl| a_N \bigl(2 (k+1)\bigr) - \biggl(\frac{3}{2 \pi} \cE(W)\biggr)^k
a_N(2)\biggr| \nonumber\\
&&\qquad=
\alpha \frac{\log N}{N} \biggl| \biggl(V_N, GV_N
\biggl((GV_N)^{2k} 1_B - \biggl(\frac{3}{2 \pi} \cE(W)\biggr)^k 1_B
\biggr)\biggr)\biggr|\\
&&\qquad = I_1 + I_2 ,
\nonumber
\end{eqnarray}
where in the last step, using the symmetry of $G$, we have set
\begin{eqnarray*}
I_1 & = & \alpha \frac{\log N}{N} \biggl| \biggl(GV_N, V_N
1_{B_k} \biggl((GV_N)^{2k} 1_B - \biggl(\frac{3}{2 \pi} \cE(W)
\biggr)^k 1_B\biggr)\biggr)\biggr| ,
\\
I_2 & = & \alpha \frac{\log N}{N} \biggl| \biggl(GV_N, V_N
1_{B\setminus B_k} \biggl((GV_N)^{2k} 1_B - \biggl(\frac{3}{2 \pi}
\cE(W)\biggr)^k 1_B\biggr)\biggr)\biggr| .
\end{eqnarray*}
We then observe that
\[
I_1 \le c \alpha\log N \|GV_N\|_\infty\|V_N\|_\infty \Delta _{k,N}
\stackrel{\mbox{\fontsize{8.36}{8.36}\selectfont{(\ref{3.10})}}}{\le}
c' \alpha\Delta_{k,N} \ulN^{\mbox{\fontsize{8.36}{8.36}\selectfont{
{(\ref{4.61})}}}} 0
\]
and that
\begin{eqnarray*}
I_2 &\stackrel{\mbox{\fontsize{8.36}{8.36}\selectfont{(\ref{4.24})}}}{\le}& \alpha\log N \|GV_N\|
_\infty\|V_N\|_\infty\times c(k) \times\frac{|B \setminus B_k|}{N}
\\
&\stackrel{\mbox{\fontsize{8.36}{8.36}\selectfont{(\ref{3.10})}}}{\le}& c'(k) \alpha
\frac{|J \setminus J_k|}{|J|}
\ulN^{\mbox{\fontsize{8.36}{8.36}\selectfont{(\ref{4.56})}}} 0.
\end{eqnarray*}
We have thus shown that
%
%
\begin{equation}\label{4.63}\quad
\lim_N \biggl| a_N \bigl(2 (k+1)\bigr) - \biggl(\frac{3}{2 \pi} \cE
(W)\biggr)^k a_N(2)\biggr| = 0\qquad \mbox{for any $k \ge1$}.
\end{equation}
Combined with~(\ref{4.48}) this completes the proof of~(\ref{4.14})
and hence of Theorem~\ref{theo4.1}.
\end{pf*}

Our next objective is to study the convergence in distribution of the
random fields $(\frac{\cL_{y,N}}{\log N})_{y \in\IZ^2}$ and $(\frac
{\cL'_{y,N}}{N u '_N})_{y \in\IZ^2}$, as $N$ goes to infinity, where
we recall the notation from~(\ref{4.1}),~(\ref{4.2}). The task is
simplified by the fact that we have already proved Theorem \ref
{theo4.2}: we only need to investigate the convergence in distribution
of these random fields at the origin. We now focus on the case where
$\Lambda= \{0\}$, and $J_0$ plays the role of $B$. We further define
%
%
\begin{equation}\label{4.64}
\wt{V}_N(x) = \frac{1}{\log N} 1_{J_0}(x),\qquad
\wt{V}'_N(x) = \frac{1}{\gamma_N} \wt{V}_N(x)\qquad\mbox{for $x \in\IZ^3$}
\end{equation}
and set
%
%
\begin{equation} \label{4.65}
\wt{\cL}_N = \frac{1}{\log N} \cL_{0,N},\qquad \wt{\cL}'_N =
\frac{1}{N u'_N} \cL'_{0,N} .
\end{equation}
Just as in~(\ref{4.7}),~(\ref{4.10}), we know by Theorem \ref
{theo2.1} that
%
%
\begin{equation}\label{4.66}
\IE[\exp\{z \wt{\cL}_N\}] = \exp\biggl\{\sum_{n \ge1} \wt{a}_N
(n) z^n\biggr\}\qquad \mbox{for $|z| < \wt{r}_N$ in $\IC$}
\end{equation}
with $r_N > 0$ and where we have set
%
%
\begin{equation}\label{4.67}
\wt{a}_N(n) = u_N (\wt{V}_N, (G\wt{V}_N)^{n-1} 1)\qquad \mbox{for $n
\ge1$}
\end{equation}
and that
%
%
\begin{equation}\label{4.68}
\IE[\exp\{z \wt{\cL}'_N\}] = \exp\biggl\{\sum_{n \ge1} \wt{a}
' _N (n) z^n\biggr\}\qquad \mbox{for $|z| < \wt{r} ' _N$ in $\IC$}
\end{equation}
with $r'_N > 0$ and
%
%
\begin{equation}\label{4.69}
\wt{a} ' _N (n) = u'_N (\wt{V}'_N, (G\wt{V}'_N)^{n-1} 1) = \frac
{1}{\alpha} \gamma_N^{1-n} \wt{a}_N(n)\qquad \mbox{for $n \ge1$}.
\end{equation}
The heart of the matter for the proof of Theorem~\ref{theo4.9} below
lies in the control of the large $N$ behavior of the sequence $\wt
{a}_N(n)$, $n \ge1$.
\begin{theorem}\label{theo4.8}
%
%
\begin{eqnarray}
\label{4.70}
&0 \le\wt{a}_N(n) \le\alpha c_7^n\qquad \mbox{for $n \ge1, N > 1$},&
\\[-2pt]
\label{4.71}
&\mbox{for any $n \ge1$}\qquad \displaystyle \lim_N \wt{a}_N(n) = \alpha
\biggl(\frac{3}{\pi}\biggr)^{n-1}.&\vspace*{-2pt}
\end{eqnarray}
\end{theorem}

We first explain how Theorems~\ref{theo4.2} and~\ref{theo4.8} enable
us to infer the convergence in law of the appropriately scaled fields
$\cL_{y,N}$, $y \in\IZ^2$, and $\cL'_{y,N}$, $y \in\IZ^2$.\vspace*{-2pt}
\begin{theorem}\label{theo4.9}
%
%
\begin{eqnarray}
\label{4.72}
\begin{tabular}{p{318pt}}
As $N$ goes to infinity, $(\frac{\cL_{y,N}}{\log N})_{y \in
\IZ^2}$ converges in distribution to a flat random
field with constant valued distributed as $R^2$ with $R$ as in
(\ref{4.18}).
\end{tabular}\hspace*{-38pt}
\\[-2pt]
\label{4.73}
\begin{tabular}{p{318pt}}
As $N$ goes to infinity, $(\frac{\cL'_{y,N}}{Nu'_N})_{y \in
\IZ^2}$ converges in distribution to a flat random
field with value $1$.
\end{tabular}\hspace*{-38pt}   \vspace*{-2pt}
\end{eqnarray}
\end{theorem}
\begin{pf*}{Proof
\textup{(assuming Theorem \ref
{theo4.8})}} A repetition\vspace*{1pt} of the arguments used in the proof of Theorem
\ref{theo4.2} shows that $\wt{\cL}_N$ converges in distribution to a
nonnegative random variable $\wt{\cL}$ with Laplace transform
%
%
\begin{equation}\label{4.74}
E[\exp\{ - \lambda\wt{\cL}\}] = \exp\biggl\{ - \frac{\alpha
\lambda}{1 + {3\lambda}/{\pi} } \biggr\}\qquad \mbox{for $\lambda
\ge0$},
\end{equation}
so that by~(\ref{1.30}), $\wt{\cL}$ is $\operatorname{BESQ}^0(\alpha, \frac{3}{2
\pi})$-distributed, that is, has the same distribution as $V^2$ in the
notation of~(\ref{4.18}).

Moreover we know from Theorem~\ref{theo4.2} that for any $y \in\IZ
^2$, when $N$ goes to infinity, $\frac{1}{\log N} (\cL_{y,N}- \cL
_{0,N})$ converges to zero in distribution, and~(\ref{4.72}) follows.

In the case of~(\ref{4.73}) we note instead that the arguments used in
the proof of Theorem~\ref{theo4.2} now show that $\wt{\cL}'_N$
converges in distribution to a nonnegative random variable with
Laplace transform $e^{-\lambda}$, $\lambda\ge0$, that is, to the
constant $1$. Since\vspace*{1pt} by Theorem~\ref{theo4.2}, $\frac{1}{N u'_N}(\cL
'_{y,N}- \cL'_{0,N})$ converges in distribution to zero for any $y \in
\IZ^2$, we obtain~(\ref{4.73}).\vspace*{-2pt}
\begin{pf*}{Proof of Theorem~\ref{theo4.8}}
We now write \mbox{$\| \cdot\|_\infty
$} for the supremum norm on $\wt{B} = J_0$, and the linear operators we
consider are restricted to functions that vanish outside $\wt{B}$. The
fact that $\wt{a}_N(n)$ is nonnegative is plain; see~(\ref{4.64}),
(\ref{4.67}). Moreover the right-hand inequality in~(\ref{4.70}) is a
direct consequence of~(\ref{3.8}). This proves~(\ref{4.70}).

We now turn to the proof of~(\ref{4.71}). For $k \ge1$, we introduce
$\wt{B}_k = \{0\} \times J_k$, with $J_k$ as in~(\ref{4.56}), so that
for any $k \ge1$, and $N \ge c(k)$, $\phi\not= \wt{B}_k \subsetneq
\wt{B}_{k-1} \subsetneq\cdots\subsetneq\wt{B}_1 \subsetneq\wt
{B}$. In a much\vadjust{\goodbreak} simpler fashion than~(\ref{4.59}) we now find that
[see~(\ref{4.39}) for notation]
%
%
\begin{equation}\label{4.75}\qquad
\bigl\|1_{\wt{B}_k} \bigl((G \wt{V}_N)^k 1_{\wt{B}} - (2 \tau_N)^k 1_{\wt
{B}}\bigr)\bigr\|_\infty\le c(k) \biggl(\|1_{\wt{B}_1} \wt{F}_N\|_\infty+
\frac{1}{\sqrt{\log N}}\biggr),
\end{equation}
where we have set
\[
\wt{F}_N = (G\wt{V}_N) 1_{\wt{B}} - 2 \tau_N 1_{\wt{B}} = G \wt
{V}_N - 2 \tau_N 1_{\wt{B}}.
\]
It already follows from the definitions of $\wt{V}_N$ and $\tau_N$ in
(\ref{4.64}),~(\ref{4.39}) that
%
%
\begin{equation}\label{4.76}
\lim_N \|1_{\wt{B}_1} \wt{F}_N\|_\infty= 0.
\end{equation}
We thus see that for $k \ge1$,
\begin{eqnarray}\label{4.77}
\lim_N \wt{a}_N(k) &=& \lim_N \alpha
\frac{\log N}{N} (\wt{V}_N, (G \wt{V}_N)^{k-1} 1) \nonumber\\[-8pt]\\[-8pt]
&=& \alpha\Bigl(2 \lim
_N \tau_N\Bigr)^{k-1} =
\alpha\biggl(\frac{3}{\pi}\biggr)^{k-1}\nonumber
\end{eqnarray}
and this proves~(\ref{4.71}).
\end{pf*}
\begin{remark}\label{rem4.10}
(1) There is an important connection between random interlacements at
level $u$ and the structure left by a random walk on a large torus
$(\IZ/ N\IZ)^d$ (here $d =3$), at a microscopic scale of order $1$
(see~\cite{Wind08}) or even at a mesoscopic scale of order
$N^{1-\varepsilon}$, with $0 < \varepsilon< 1$ (see
\cite{TeixWind10}) when the walk runs for times of order $uN^d$. This
naturally raises the question whether the above limiting results might
also be relevant for the field of occupation times left close to the
origin, by continuous time simple random walk with uniform starting
point on a large two-dimensional torus $(\IZ/N\IZ)^2$, at times of
order $\alpha N^2 \log N$ $(\mbox{$=$} u_N N^3)$ or at much larger times
$u'_N N^3$. Let us incidentally point out that the time scale $\alpha
N^2 \log N$ is much smaller than the cover time of the torus which has
order $\frac{4}{\pi} N^2 (\log N)^2$; see~\cite{DembPereRoseZeit04}.

(2) We can consider the discrete occupation times $\ell_{x,u}$, $x \in
\IZ^d$, $u \ge0$ [see Remark~\ref{rem2.4}(5)] and define
${\mathfrak L}_{y,N}$ and ${\mathfrak L}'_{y,N}$, for $y \in\IZ^2$,
$N > 1$, as in~(\ref{4.2}), simply replacing $L_{x,u}$ by $\ell_{x,u}$.

Theorems~\ref{theo4.2} and~\ref{theo4.9} enable us to show that
$(\frac{{\mathfrak L}'_{y,N} - {\mathfrak L}'_{0,N}}{\sqrt{N
u'_N}})_{y \in\IZ^2}$ converges in distribution to a centered
Gaussian field which vanishes at the origin. However, this limiting
field is different from $(\psi_y)_{y \in\IZ^2}$ in~(\ref{4.20}).
\end{remark}

The heart of the matter lies in the fact that for any $W(\cdot)$ as in
(\ref{4.3}), when $N$ tends to infinity,
%
%
\begin{equation}\label{4.78}
\begin{tabular}{p{318pt}}
$\sum_{y \in\IZ^2} W(y) \frac{{\mathfrak L}'_{y,N}}{\sqrt
{N u'_N}}$ converges in distribution to a Gaussian variable
with variance $\cE(W) - \sum_{y \in\IZ^2} W(y)^2$
\end{tabular}\hspace*{-36pt}
\end{equation}
[in the case of $\cL'_{y,N}$ the limiting variance instead equals $\cE(W)$].\vadjust{\goodbreak}

Indeed it follows from~(\ref{2.27}) that for real $z$ and $N > 1$,
\begin{eqnarray}\label{4.79}
&&\IE\biggl[\exp\biggl\{ z \sum_{y \in\IZ^2} \frac{W(y)}{\sqrt
{Nu'_N}} {\mathfrak L}'_{y,N}\biggr\}\biggr]\nonumber\\[-8pt]\\[-8pt]
&&\qquad = \IE\biggl[\exp\biggl\{
\sum_{y \in\IZ^2} \bigl(1 - e^{-z {W(y)}/{\sqrt{N
u'_N}}}\bigr) \cL'_{y,N}\biggr\}\biggr],\nonumber
\end{eqnarray}
and for $|z| < r$ and $N \ge c$, we can use Taylor's expansion and write
\[
1 - e^{-z{W(y)}/{\sqrt{Nu'_N}}} = z
\frac{W(y)}{\sqrt{N u'_N}} - \fr z^2 \frac{W(y)^2}{Nu'_N}
\bigl(1 + \varepsilon_y(z,N)\bigr),
\]
where $|\varepsilon_y (z,N)| \le\frac{1}{2}$ and $\lim_N
\varepsilon_y(z,N) = 0$, for each $y \in\IZ^2$, $|z| < r$.

Inserting the above identity in~(\ref{4.79}) shows that for $|z| < r$
and $N \ge c$,
the left-hand side of~(\ref{4.79}) equals
%
%
\begin{equation}\label{4.80}
\IE\biggl[\exp\biggl\{ z \sum_{y \in\IZ^2} \frac{W(y)}{\sqrt
{Nu'_N}} \cL'_{y,N} - \fr z^2 \sum_{y \in\IZ^2} \frac
{W(y)^2}{Nu'_N} \bigl(1 + \varepsilon_y(z,N)\bigr)
\cL'_{y,N}\biggr\}\biggr].\hspace*{-28pt}
\end{equation}
By the end of the proof of Theorem~\ref{theo4.2} we know that for $|z|
< c$,
\[
\lim_N \IE\biggl[\exp\biggl\{ z \sum_{y \in\IZ^2}
\frac{W(y)}{\sqrt{Nu'_N}} \cL'_{y,N}\biggr\}\biggr] = \exp\biggl\{
\frac{z^2}{2} \cE(W)\biggr\}.
\]
A straightforward uniform integrability argument combined with (\ref
{4.73}) and (\ref{4.80}) shows that for real $z$ with $|z| < c$,
%
%
\begin{equation}\label{4.81}
\lim_N \IE\biggl[\exp\biggl\{ z \sum_{y \in\IZ^2}
\frac{W(y)}{\sqrt{Nu'_N}} {\mathfrak L}'_{y,N}\biggr\}\biggr] = \exp
\biggl\{\frac{z^2}{2} \cE(W) - \frac{z^2}{2} \sum_{y \in\IZ
^2} W(y)^2\biggr\}.\hspace*{-28pt}
\end{equation}
Similar arguments as in the proof of Theorem~\ref{theo4.2} now yield
(\ref{4.78}).

Note incidentally that Theorems~\ref{theo4.2} and~\ref{theo4.9} are
not quite sufficient to study the limit in law of $(\frac{{\mathfrak
L}_{y,N} - {\mathfrak L}_{0,N}}{\sqrt{\log N}} )_{y \in\IZ^2}$. As
shown by the above proof [see in particular~(\ref{4.80})], to handle
this case we would in essence need a limiting result for the joint law
of the two random fields that appear in~(\ref{4.17}) and~(\ref{4.72}).
\end{pf*}

\section{Occupation times at high level $u$}\label{sec5}

In this section we relate occupation times at a high level $u$ of the
random interlacements with the $d$-dimensional Gaussian free field. The
limit $u \rightarrow\infty$ bypasses the obstructions present when one
considers a fixed level $u$; see Remark~\ref{rem2.4}(2). Our main
result appears in Theorem~\ref{theo5.1}. It has a similar flavor to
(\ref{4.20}) of Theorem~\ref{theo4.2} and~(\ref{4.73}) of Theorem
\ref{theo4.9}. Moreover it can rather straightforwardly be extended to
the case of random interlacements on transient weighted graphs; see
Remark~\ref{rem5.2}. However, we keep the set-up of $\IZ^d$, $d \ge
3$, for the main body of this section, not to overburden notation.

We consider on an auxiliary probability space
%
%
\begin{equation}\label{5.1}
\begin{tabular}{p{310pt}}
$\varphi_x, x \in\IZ^d$, a centered Gaussian field with
covariance function
$E[\varphi_x \varphi_{x'}] = g(x' - x) + g(0) - g(x) - g(x')$,
for $x,x' \in\IZ^d$.
\end{tabular}\hspace*{-32pt}
\end{equation}
This field has the same distribution as the field $(\gamma_x - \gamma
_0)_{x \in\IZ^d}$ of increments at the origin of the $d$-dimensional
Gaussian free field, $(\gamma_x)_{x \in\IZ^d}$, that is, the
centered Gaussian field with covariance function $E[\gamma_x \gamma
_{x'}] = g(x,x')$, for $x,x' \in\IZ^d$.

We can now state the main result of this section.
\begin{theorem}\label{theo5.1}
As $u \rightarrow\infty$,
%
%
\begin{equation} \label{5.2}
\begin{tabular}{p{310pt}}
$(\frac{L_{x,u} - L_{x,0}}{\sqrt{2 u}})_{x \in
\IZ^d}$ converges in distribution to the Gaussian
random field $(\varphi_x)_{x \in\IZ^d}$ in
(\ref{5.1})
\end{tabular}\hspace*{-32pt}
\end{equation}
and
%
%
\begin{equation} \label{5.3}
\begin{tabular}{p{310pt}}
$\biggl(\dfrac{1}{u} L_{x,u}\biggr)_{x \in\IZ^d}$ converges in
distribution to the constant field equal to $1$.
\end{tabular}\hspace*{-32pt}
\end{equation}
\end{theorem}
\begin{pf}
We follow the same strategy as in the previous section, the problem is,
however, much simpler now. It clearly suffices to prove~(\ref{5.2})
and~(\ref{5.3}) with $u$ replaced by a sequence $u_N$ such that
%
%
\begin{equation}\label{5.4}
u_N \ge1\qquad \mbox{for } N \ge1\quad \mbox{and}\quad \lim_N u_N = \infty.
\end{equation}
We thus consider a function $V$ on $\IZ^d$ such that
%
%
\begin{equation}\label{5.5}
V\mbox{ is finitely supported}\quad\mbox{and}\quad\sum_{x \in\IZ^d} V(x) = 0.
\end{equation}
We define
%
%
\begin{equation}\label{5.6}
L_N = \sum_{x \in\IZ^d} \frac{1}{\sqrt{2 u_N}} V(x)
L_{x,u_N}\qquad \mbox{for $N \ge1$}.
\end{equation}
It follows from Theorem~\ref{theo2.1} and Remark~\ref{rem2.2} that
for some fixed $r > 0$,
%
%
\begin{equation}\label{5.7}
\IE[\exp\{z L_N\}] = \exp\biggl\{\sum_{n \ge1} c_N(n) z^n\biggr\}
\qquad\mbox{for $|z| < r$ in $\IC$},
\end{equation}
where we have set for $n, N \ge1$,
%
%
\begin{equation}\label{5.8}
c_N(n) = 2^{-n/2} u^{1-{n/2}}_N (V,(GV)^{n-1} 1).
\end{equation}
In view of~(\ref{5.5}) we have
%
%
\begin{equation}\label{5.9}
c_N(1) = 0.
\end{equation}
By~(\ref{5.4}),~(\ref{5.8}) it is also plain that for $n \ge2$,
%
%
\begin{equation}\label{5.10}
|c_N(n)| \le c(V)^n
\end{equation}
with $c(V)$ a positive constant depending on $V$ and $d$, by our
convention. Moreover we find that
%
%
\begin{equation}
\label{5.11}
\lim_N c_N(n) = 0\qquad \mbox{for } n > 2
\end{equation}
and
%
\begin{equation}
\label{5.12}
c_N(2) = \fr(V, GV) \stackrel{\mbox{\fontsize{8.36}{8.36}\selectfont{(\ref{5.1}),~(\ref{5.5})}}}{=} \fr
E\biggl[\biggl(\sum_x V(x) \varphi_x\biggr)^2\biggr].
\end{equation}
The same arguments as in the proof of Theorem~\ref{theo4.2} show that
%
%
\begin{eqnarray}\label{5.13}
\begin{tabular}{p{310pt}}
$L_N$ converges in distribution to a centered Gaussian variable
with variance
$E[(\sum_x V(x) \varphi_x)^2]$.
\end{tabular}\hspace*{-36pt}
\end{eqnarray}
Since $V$ in~(\ref{5.5}) and $u_N$ in~(\ref{5.4}) are arbitrary,
claim~(\ref{5.2}) follows. We then turn to the proof of~(\ref{5.3}).
It follows by~(\ref{2.21}) that for $\lambda\ge0$,
\begin{eqnarray}\label{5.14}
\IE\biggl[\exp\biggl\{- \frac{\lambda}{u_N} L_{0,u_N}\biggr\}\biggr]
&=& \exp\biggl\{ - \frac{\lambda}{1 + g(0) ({\lambda
}/{u_N})}\biggr\} \nonumber\\[-8pt]\\[-8pt]
&\longrightarrow& e^{-\lambda}\qquad\mbox{as } N \rightarrow
\infty.\nonumber
\end{eqnarray}
This shows that $\frac{1}{u_N} L_{0,u_N}$ converges in distribution
to the constant $1$ as $N$ goes to infinity. Since\vspace*{2pt} due to~(\ref{5.2}),
for any $x \in\IZ^d$, $\frac{1}{u_N} (L_{x,u_N} - L_{0,u_N})$ tends
to zero in distribution, as $N$ goes to infinity, our claim follows.
\end{pf}
\begin{remark}\label{rem5.2}
(1) The results of the present section can straightforwardly be
extended to the set-up of continuous time random interlacements on a
transient weighted graph $E$, as we now explain. We keep the same
notation and assumptions as in Remark~\ref{rem2.4}(4). We introduce a base
point $x_0 \in E$. In place of~(\ref{5.1}) we consider
%
%
\begin{equation}\label{5.15}
\begin{tabular}{p{305pt}}
$\varphi_x, x \in E$, a centered Gaussian field with covariance
function
$E[\varphi_x \varphi_{x'}] = g(x,x') + g(x_0,x_0) - g(x_0,x) -
g(x_0,x'), x,x' \in E$,
\end{tabular}\hspace*{-32pt}
\end{equation}
where $g(\cdot,\cdot)$ now stands for the Green density.

This field has the same distribution as the field of increments
$(\gamma_x - \gamma_{x_0})_{x \in E}$, of the Gaussian free field
$(\gamma_x)_{x \in E}$ attached to the transient weighted graph, that
is, the centered Gaussian field with covariance function $E[\gamma_x
\gamma_{x'}] = g(x,x')$, $x,x' \in E$.

With the help of~(\ref{2.25}),~(\ref{2.26}), the arguments employed
in the proof of Theorem~\ref{theo5.1} now show that as $u \rightarrow
\infty$,
%
%
\begin{equation}\label{5.16}
\begin{tabular}{p{310pt}}
$(\frac{L_{x,u} - L_{x,0}}{\sqrt{2 u}})_{x \in
E}$ converges in distribution to the Gaussian
random field $(\varphi_x)_{x \in E}$
\end{tabular}\hspace*{-32pt}
\end{equation}
and that
%
%
\begin{equation}\label{5.17}
\biggl(\frac{1}{u} L_{x,u}\biggr)_{x \in E} \!\mbox{ converges in
distribution to the constant field equal to $1$}.\hspace*{-35pt}
\end{equation}

(2) In the case of the discrete occupation times $\ell_{x,u}$, $x \in
\IZ^d$, $u \ge0$, the same arguments used in Remark~\ref{4.10}(2)
show that when $u \rightarrow\infty$,
%
%
\begin{equation}\label{5.18}
\biggl(\frac{\ell_{x,u} - \ell_{x,0}}{\sqrt{2 u}}\biggr)_{x \in
\IZ^d} \mbox{ converges in distribution to $(\nu_x)_{x \in\IZ^d}$},
\end{equation}
where $(\nu_x)_{x \in\IZ^d}$ is the centered Gaussian field
vanishing at the origin such that for any $V$ as in~(\ref{5.5}), one
has in the notation of~(\ref{5.1})
%
%
\begin{equation}\label{5.19}
E\biggl[\biggl(\sum_{x \in\IZ^d} V(x) \nu_x\biggr)^2\biggr] + \fr
\sum_{x \in\IZ^d} V(x)^2 = E\biggl[\biggl(\sum_{x \in\IZ^d} V(x)
\varphi_x\biggr)^2\biggr].
\end{equation}
Moreover looking at the Laplace functional, one sees with the help of
(\ref{2.27}) and~(\ref{5.3}) that for $u \rightarrow\infty$,
%
%
\begin{equation}\label{5.20}
\biggl(\frac{1}{u} \ell_{x,u}\biggr)_{x \in\IZ^d} \!\mbox
{ converges in distribution to the constant field equal to $1$}.\hspace*{-35pt}
\end{equation}
\end{remark}


%

%
\printaddresses

\end{document}